\documentclass[11pt]{amsart}

\usepackage{amsmath, amssymb, amsthm, latexsym, amscd}

\newtheorem{theorem}[subsection]{Theorem}
\newtheorem{prop}[subsection]{Proposition}
\newtheorem{lemma}[subsection]{Lemma}
\newtheorem{corollary}[subsection]{Corollary}

\theoremstyle{definition}
\newtheorem{define}[subsection]{Definition}
\newtheorem{notation}[subsection]{Notation}

\newtheorem{remark}[subsection]{Remark}

\newcommand{\mf}[1]{\mathfrak{{#1}}}

\newcommand{\C} {\mathbb{C}}

\newcommand{\V}{\mathbb{V}}

\newcommand{\A}{\mathcal{A}}

\title[VOA structure on modified regular representations]{vertex operator algebras associated to modified regular representations of affine Lie algebras}
\author{Minxian Zhu}
\address{Department of Mathematics, Yale University, New Haven, CT 06520}
\email{minxian.zhu@yale.edu}

\begin{document}
\maketitle

\begin{abstract}

Let $G$ be a simple complex Lie group with Lie algebra $\mf g$ 
and let $\hat{\mathfrak{g}}$ be the affine Lie algebra. 
We use intertwining operators and Knizhnik-Zamolodchikov equations 
to construct a family of $\mathbb{N}$-graded vertex operator algebras associated to $\mf g$. 
They are $\hat{\mf g}\oplus \hat{\mf g}$-modules of dual levels $k, \bar k \notin \mathbb Q$ 
in the sense that $k + \bar k = -2 h^\vee$ where $h^{\vee}$ is the dual Coxeter number of $\mf g$. 
Its conformal weight $0$ component is the algebra of regular functions on $G$. 
This family of vertex operator algebras were previously studied by Arkhipov-Gaitsgory and Gorbounov-Malikov-Schechtman from different points of view. 
We show that the vertex envelope of the vertex algebroid associated to $G$ and level $k$
is isomorphic to the vertex operator algebra we constructed above when $k$ is irrational. 
The case of integral central charges is also discussed.

\end{abstract}

\section{introduction}

Let $\mf{g}$ be a simple complex Lie algebra, $\hat{\mf{g}}$ be the affine Lie algebra. 
The aim of this paper is to study a one-parameter family of vertex operator algebras associated to $\mf{g}$, which have the structure of  $\hat{\mf{g}} \oplus  \hat{\mf{g}}$-representations with dual central charges as in the works of [AG], [GMS1,2] and [FS].

In [FS], it is constructed explicitly for $\widehat{sl_{2}}$ using Fock space realizations and it decomposes into summands corresponding to the dominant integral weights of $\mf{g}$ for generic central charges. More specifically one has the following 
$$\V =\bigoplus_{\lambda \in P^{+}}V_{\lambda,k} \otimes V_{\lambda^{*}, \bar{k}}$$
for $k, \bar{k}\notin \mathbb{Q}$ and dual in the sense that $k+\bar{k} = -2h^{\vee}$, where $h^{\vee}$ is the dual Coxeter number of $\mf{g}$. 
Here $V_{\lambda, k} = U(\hat{\mf{g}})\otimes_{U(\mf{g}\otimes  \mathbb{C}[t]\oplus \mathbb{C} \underline c)} V_{\lambda}$ is the so called Weyl module induced from $V_{\lambda}$, the irreducible representation of $\mf{g}$ with highest weight $\lambda$, by letting $\mf{g}\otimes t\mathbb{C}[t]$ act by $0$ and the central element $\underline c$ act as multiplication by the scalar $k$ on $V_{\lambda}$.  
The module induced from $V_{\lambda}^{*}$ of level $\bar{k}$ is denoted by $V_{\lambda^{*}, \bar{k}}$. 

$\V$ is naturally $\mathbb{N}$-graded and the top level $\V_{0} = \oplus_{\lambda \in P^{+}} V_{\lambda}\otimes V_{\lambda}^{*}$ can be identified with $\mathcal{R}(G)$, the space of regular functions on the simple complex Lie group $G$ associated with $\mf{g}$. 
According to [FZ] the vacuum module $V_{0,k}$ (resp. $V_{0,\bar{k}}$) carries the structure of a vertex operator algebra(VOA) and other modules $V_{\lambda,k}$ (resp. $V_{\lambda^{*}, \bar{k}}$) become its representations. 
In fact the zero block $V_{0, k}\otimes V_{0,\bar{k}}$ is a vertex subalgebra of $\V$ and it generates the two copies of $\hat{\mf{g}}$-actions on $\V$. 
For generic values of $k$ and $\bar{k}$, all the modules involved are irreducible and the space of intertwining operators of the type $\left( \begin{array}{ccc} \quad & V_{\nu,k} & \quad \\  V_{\lambda,k} & \quad & V_{\mu,k} \end{array} \right)$ (resp. $\left( \begin{array}{ccc} \quad & V_{\nu^{*},\bar{k}} & \quad \\  V_{\lambda^{*},\bar{k}} & \quad & V_{\mu^{*},\bar{k}} \end{array} \right)$) is isomorphic to $\textrm{Hom}_{\mf{g}}(V_{\lambda}\otimes V_{\mu}, V_{\nu})$ (resp. $\textrm{Hom}_{\mf{g}}(V_{\lambda}^{*}\otimes V_{\mu}^{*}, V_{\nu}^{*})$). 
In Section 2 we construct the vertex operators explicitly by pairing all these intertwining operators together in an appropriate manner (cf. [St]). 
More specifically the structural coefficients are the same as those to describe the multiplication in $\mathcal{R}(G)$ in terms of intertwining operators for the left and right $\mf{g}$-actions.  
We show the locality of vertex operators by using some knowledge of Knizhnik-Zamolodchikov equations and therefore prove the following.

\newcounter{foo}
\newtheorem{localtheorem}[foo]{Theorem}

\begin{localtheorem}
Let $\mf{g}$ be a simple complex Lie algebra and $\varkappa \notin \mathbb{Q}$. Set $k=\varkappa - h^{\vee}$, $\bar{k} = - \varkappa - h^{\vee}$, then 
$$\V= \bigoplus_{\lambda \in P^{+}}V_{\lambda,k} \otimes V_{\lambda^{*},\bar{k}}$$ 
is a vertex operator algebra of rank $2\,\textrm{dim}(\mf{g})$. 
\end{localtheorem}

The second half of the paper is based on the work of Gorbounov, Malikov and Schechtman about the chiral algebras of differential operators (cf. [GMS1], [GMS2]). 
They introduced the notion of a vertex algebroid which basically captures the structure inherited by the top two levels of an $\mathbb{N}$-graded vertex algebra. 
One can assign a vertex algebroid to a vertex algebra by truncation and conversely 
it is shown in [GMS1] that one can construct a vertex algebra from a vertex algebroid and the reconstruction functor is left adjoint to the truncation functor. 
In fact an $\mathbb{N}$-graded conformal algebra is obtained first from the vertex algebroid, then the vertex envelope of the conformal algebra is constructed which already carries the structure of a vertex algebra, and finally the quotient of the vertex envelope by a vertex ideal gives the enveloping algebra of the vertex algebroid. 
An important example in [GMS2] is the vertex algebroid $\mathcal{A}_{\mf{g},k}$ associated to the group $G$ and a fixed level $k$.
Its zeroth level is the algebra of regular functions on $G$ and the first level is the direct sum of vector fields $T$ and one forms $\Omega$ (as a vector space). 
Applying the construction functor we get a vertex algebra denoted by $U\mathcal{A}_{\mf{g}, k}$. Besides the obvious embedding $V_{0,k}\hookrightarrow U\mathcal{A}_{\mf{g}, k}$, 
the vertex algebra $V_{0,\bar{k}}$ of the dual level can also be imbedded into $U\mathcal{A}_{\mf{g}, k}$ and their images commute in the appropriate sense (cf. [AG], [GMS2]).  
In particular $U\mathcal{A}_{\mf{g}, k}$ is a $\hat{\mf{g}} \oplus \hat{\mf{g}}$-module of dual levels $k$ and $\bar{k}$.

In Section 3 we further analyze the structure of $U\mathcal{A}_{\mf{g}, k}$ and describe it using generating fields and operator product expansions (OPE).  
Set $A=\mathcal{R}(G)$. 
Let $\{\tau_{i}\}$ be an orthonormal basis of $\mf{g}$ with respect to the normalized invariant bilinear form on $\mf{g}$ which act on $A$ as left invariant vector fields and let $\{\omega_{i}\}$ be the left invariant one forms dual to $\{\tau_{i}\}$. 
Let $W$ be the vertex algebra generated by the quantum fields
$a(z) = \sum a_{(n)} z^{-n-1}$, $a\in A$ of conformal weight $0$ and 
$\tau_i (z) = \sum {\tau_i}_{(n)} z^{-n-1}$, 
$\omega_i (z) = \sum {\omega_i}_{(n)} z^{-n-1}$ of conformal weight $1$ with the OPE: 
$$\tau_{i}(z_{1}) a(z_{2}) \sim \frac{\tau_{i} a (z_{2})}{z_{1}-z_{2}}$$  
$$\tau_{i}(z_{1})\tau_{j}(z_{2})\sim \frac{C_{ijs}\tau_{s}(z_{2})}{z_{1}-z_{2}} + \frac{k\delta_{ij}}{(z_{1}-z_{2})^{2}}$$ 
$$\tau_{i}(z_{1})\omega_{j}(z_{2})\sim \frac{C_{ijs}\omega_{s}(z_{2})}{z_{1}-z_{2}} + \frac{\delta_{ij}}{(z_{1}-z_{2})^{2}}$$
$$a(z_{1})b(z_{2})\sim 0 \qquad  a(z_{1})\omega_{j}(z_{2})\sim 0 \qquad \omega_{i}(z_{1})\omega_{j}(z_{2})\sim 0$$
where $C_{ijk}$'s are the structure coefficients of $\mf{g}$ with respect to the basis $\{\tau_{i}\}$.  
Let $I$ be the vertex $\partial$-ideal generated by 
$a_{(-1)} b_{(-1)} \mathbf 1 - (a b )_{(-1)} \mathbf 1$,  
$a_{(-2)} \mathbf 1 - \sum_i (\tau_i a)_{(-1)} {\omega_i}_{(-1)} \mathbf 1$ 
where $a, b\in A$ and $a b$ denotes the multiplication of functions, 
then $W/I \cong U\mathcal{A}_{\mf{g}, k}$. 
Certainly we can replace the left invariant vector fields or one forms by the right invariant ones and make corresponding changes of the OPE and generators of the vertex ideal $I$ to formulate three other descriptions. 
Furthermore there is a Virasoro element $\varpi\in {U\mathcal{A}_{\mf{g}, k}}_{2}$ 
satisfying $\varpi_{(0)}\varpi = \partial \varpi$, $\varpi_{(1)} \varpi = 2 \varpi$, $ \varpi_{(2)}\varpi =0$
and $\varpi_{(3)} \varpi = \textrm{dim}\mf{g} $, 
hence $U\mathcal{A}_{\mf{g}, k}$ is in fact a vertex operator algebra of rank $2\,\text{dim} \mf{g}$. 
We show that the Zhu's algebra for $U\mathcal{A}_{\mf{g}, k}$ is isomorphic to the algebra of differential operators on the Lie group $G$. 
We also introduce a nondegenerate symmetric bilinear form $(, )$ on $U\mathcal{A}_{\mf{g}, k}$ compatible with the VOA structure in the sense of [FHL]. 
When $k$ is irrational, $U\mathcal{A}_{\mf{g}, k}$ coincides with the vertex operator algebra $\V$ constructed in Section 2. 

\begin{localtheorem}
Let $\mf{g}$ be a simple complex Lie algebra and $k\notin \mathbb{Q}$. 
Let $\mathcal{A}_{\mf{g}, k}$ be the associated vertex algebroid, 
then as a $\hat{\mf{g}}\oplus \hat{\mf{g}}$-module 
$$U\mathcal{A}_{\mf{g}, k} \cong \bigoplus_{\lambda\in P^{+}} V_{\lambda, k}\otimes V_{\lambda^{*}, \bar{k}}$$
where $\bar{k}=-2h^{\vee}-k$, moreover $U\mathcal{A}_{\mf{g}, k} \cong \V$ as  vertex operator algebras.
\end{localtheorem}

In Section 4 we discuss the case of integral central charges. 
The construction of the enveloping algebra $U\mathcal{A}_{\mf{g}, k}$ in [GMS1, 2] works for any $k$ (while the Fock space realization in [FS] exists for any nonzero $\varkappa$, i.e. excluding the critical level $-h^\vee$). 
In particular the "size" of the vertex operator algebra $U\mathcal{A}_{\mf{g}, k}$ does not change, 
but as a $\hat{\mf{g}}\oplus \hat{\mf{g}}$-representation 
it no longer decomposes nicely into summands corresponding to dominant weights, 
instead it has linkings determined by the shifted action of the affine Weyl group. 
In particular $U\mathcal{A}_{\mf{g}, k}$ is not generated from the top level $A = \mathcal{R}(G)$ as a $\hat{\mf{g}} \oplus \hat{\mf{g}}$-module. 
We can get a glimpse of this phenomenon based on the computations done in Section 3. 
Let $\mathcal{O}_{\varkappa}$ be the category of $\mf{g}$-integrable $\hat{\mf{g}}$-modules of finite length of central charge $\varkappa -h^{\vee}$ and 
$\mathcal{O}_{q}$ be the category of finite dimensional representations of the quantum group $U_{q}(\mf{g})$, 
where $q=\textrm{exp} (i \pi /m\varkappa)$
and $m$ is the ratio of the length of a long and a short root of $\mf{g}$, squared. 
D. Kazhdan and G. Lusztig proved the equivalence of tensor categories $\mathcal{O}_{\varkappa}$ and $\mathcal{O}_{q}$ for $\varkappa \notin  \mathbb{Q}_{\geq 0}$ (cf. [KL1-4]). 
Meanwhile as an easy consequence of a strong semi-infinite duality in [A], 
S.M.Arkhipov proved that the category of all modules with a Weyl filtration in level $k$ is contravariantly equivalent to the analogous category in dual level $-2h^{\vee}-k$. 
Under this equivalence projective objects in positive level are transformed into tilting modules in negative level (cf. [So]). 
These observations lead us to the conjecture that the $\hat{\mf{g}}\oplus \hat{\mf{g}}$-module structure of the vertex operator algebra $U\mathcal{A}_{\mf{g}, k}$ with integral central charges (not critical) is the same as the bimodule structure of the regular representation of the corresponding quantum group at roots of unity (cf. [Z]) under equivalence of categories.  

My deepest gratitude goes to my advisor Igor Frenkel for his guidance, generosity and patience. I'd also like to thank  Konstantin Styrkas for many helpful discussions. 

\section{Construction from representations}

Let us first recall what a vertex operator algebra is, its representations and the intertwining operators between them. See [B, FB, FHL, FLM, K, LL] for more details and variations in definitions. 

\begin{define}
A vertex algebra is an $ \mathbb{N} $-graded vector space $V=\bigoplus_{n=0}^{\infty}V_n$ with a vacuum vector $\mathbf{1}\in V_0$, a linear operator $\partial$ of degree one and a linear map  
$$Y(\cdot, z): V\to \text{End}V[[z^{\pm1}]];     \qquad           
v \mapsto Y(v, z) = \sum_{n\in \mathbb{Z}} v_{(n)} z^{-n-1}$$ 
 such that for any $u, v \in V$, 
 $$\text{deg} \,v_{(n)} = -n-1+\text{deg}\, v  \quad \text{for homogeneous   }  v,$$ 
$$Y(\mathbf{1}, z) =\text{Id  }, \quad Y(v, z) \mathbf{1} \in V[[z]] \,\,\text{and} \,\, \lim_{z\to 0}Y(v, z) \mathbf{1} = v, $$
$$[\partial, Y(v, z)] = \partial_{z} Y(v, z), $$
and the Jacobi identity holds
$$z_{0}^{-1} \delta \Big( \frac{z_{1}-z_{2}}{z_{0}} \Big) Y(u, z_{1}) Y(v, z_{2}) - z_{0}^{-1} \delta\Big( \frac{ z_{2}- z_{1}}{-z_{0}} \Big) Y(v, z_{2}) Y(u, z_{1})$$
$$ = z_{2}^{-1} \delta \Big( \frac{z_{1}- z_{0}}{z_{2}} \Big)  Y( Y(u,z_{0}) v, z_{2} ).$$
Note that the Jacobi identity can be replaced by the locality condition (cf. [FB], [LL]) that for any $u,v\in V$ there exists an $N(u,v) \in \mathbb{N} $ such that 
$$(z_{1}-z_{2})^{N(u,v)} [Y(u, z_{1}), Y(v, z_{2})]=0. $$

If in addition there exists a \emph{Virasoro element} $\omega \in V_{2}$ with $$Y(\omega, z)=\sum_{n\in \mathbb{Z}} L_{n} z^{-n-2}$$ such that 
$$L_{0} | _{V_{n}} = n \text{ Id}, \quad L_{-1} = \partial, $$
$$[L_{m}, L_{n}]= (m-n)L_{m+n}+\delta_{m+n,0}\frac{m^3-m}{12}c \,\text{ Id}, $$  
V is called a vertex operator algebra of rank $c$. 
\end{define}

From the Jacobi identity we have Borcherds' commutator formula and iterate formula: 
$$[u_{(m)}, v_{(n)}]=\sum_{k\geq 0} {m \choose k} (u_{(k)}v)_{(m+n-k)}$$
$$(u_{(m)}v)_{(n)} w = \sum_{i\geq 0}(-1)^{i} {m \choose i} (u_{(m-i)}v_{(n+i)} w-(-1)^{m}v_{(m+n-i)}u_{(i)} w)$$
for any $u, v, w \in V$, $m, n \in \mathbb{Z}$.

A vertex algebra $V$ is called abelian if $u_{(n)}v=0$ for all $u,v\in V$, $n\geq 0$. 
An abelian vertex algebra is an $\mathbb{N}$-graded commutative associative algebra with unity 
$\mathbf{1}$ and a derivation $\partial$ of degree one.

\begin{define}
Given a vertex operator algebra $V$ of rank $c$, 
a representation of $V$ is an $\mathbb{N}$-graded vector space $M=\bigoplus_{n=0}^{\infty} M_{n}$ equipped with a linear map 
$Y_M (\cdot, z): V \to  \text{End}M[[z,z^{-1}]]; v\mapsto Y_{M}(v ,z)=\sum_{n\in \mathbb{Z}} v_{(n)} z^{-n-1}$  such that
$$ \text{deg}\, v_{(n)} = -n-1+\text{deg}\, v  \quad \text{for homogeneous   }  v,$$
$$ Y_{M}(\mathbf{1}, z) = \text {Id},  \quad Y_{M}( \partial v, z) = \partial _{z} Y_{M} (v, z),$$
$$z_{0}^{-1} \delta \Big( \frac{z_{1}-z_{2}}{z_{0}} \Big) Y_{M} (u, z_{1}) Y_{M} (v, z_{2}) - z_{0}^{-1} \delta\Big( \frac{ z_{2}- z_{1}}{-z_{0}} \Big) Y_{M} (v, z_{2}) Y_{M} (u, z_{1}) $$
$$ = z_{2}^{-1} \delta \Big( \frac{z_{1}- z_{0}}{z_{2}} \Big)  Y_{M} ( Y(u,z_{0}) v, z_{2} )$$
for any $u, v\in V$ and 
$$[L_{m}, L_{n}]=(m-n)L_{m+n} + \delta_{m+n,0} \frac{m^3-m}{12}c \, \text{Id}$$ 
where
$Y_{M} (\omega, z) =\sum _{n\in \mathbb{Z}} L_{n} z^{-n-2}$. 
\end{define}

\begin{define}
Given a vertex operator algebra $V$ and $V$-representations $M^{1}, M^{2}$ and $ M^{3}$, 
an intertwining operator of type $\left( \begin{array}{ccc} \quad & M^{3} & \quad \\  M^{1} & \quad & M^{2} \end{array} \right)$ is a linear map
$$\Phi(\cdot, z): M^{1} \to \text{End} (M^{2}, M^{3})\{z\}$$ such that
$$\Phi(L_{-1} v, z)=\partial_z \Phi (v, z),$$
and 
$$z_{0}^{-1} \delta \Big( \frac{z_{1}-z_{2}}{z_{0}} \Big) Y_{M^3} (u, z_{1}) \Phi (v, z_{2}) - z_{0}^{-1} \delta\Big( \frac{ z_{2}- z_{1}}{-z_{0}} \Big) \Phi (v, z_{2}) Y_{M^2} (u, z_{1}) $$
$$ = z_{2}^{-1} \delta \Big( \frac{z_{1}- z_{0}}{z_{2}} \Big)  \Phi ( Y_{M^1}(u, z_{0}) v, z_{2} ).$$
for any $u\in V, v \in M^{1}$. 
We denote the vector space of all intertwining operators of the above type by $\mathcal{V}_{M_{1} M_{2}}^{M_{3}}$. 
\end{define}

Let $G$ be any simple complex Lie group with Lie algebra $\mathfrak{g}$. 
Let $P^{+}$, $W$, $\omega_{0}$, $h^{\vee}$, $\rho$ denote the dominant integral weights, Weyl group, element of maximal length in the Weyl group, dual Coxeter number and half the sum of all positive roots respectively. 
Let $(,)$ be the normalized invariant symmetric bilinear form on $\mathfrak{g}$ which equals $\frac{1}{2h^{\vee}} (, )_{\textit{Killing}}$.

Let $\hat{\mathfrak{g}}=\mathfrak{g}\otimes \mathbb{C}[t, t^{-1}] \oplus \mathbb{C} \underline c$ be the affine Lie algebra which is naturally $\mathbb{Z}$-graded: $\text{deg} (\mathfrak{g} \otimes t^{n}) = -n, \text{deg} (\underline c) = 0$. 
Define $\hat{\mathfrak{g}}_{+}=\mathfrak{g} \otimes t\mathbb{C}[t]$,  
$\hat{\mathfrak{g}}_{-}=\mathfrak{g} \otimes t^{-1}\mathbb{C}[t^{-1}]$,  
$\hat{ \mathfrak{g}}_{0}=\mathfrak{g} \oplus \mathbb{C} \underline c$ and 
$ \hat{\mathfrak{g}}_{\geq 0}=\hat{\mathfrak{g}}_{+} \oplus \hat{\mathfrak{g}}_{0}$.  
For $x\in \mf{g}$ denote $x \otimes t^{n}$ by $x(n)$, 
then the commutators in $\hat{\mf g}$ are given by 
$$[x(m),y(n)]=[x,y](m+n)+m \delta_{m+n,0}(x,y) \underline c, $$
$$[x(n), \underline c] = [ \underline c , \underline c] = 0. $$ 

Let $V_{\lambda}$ denote the (finite dimensional) irreducible representation of $\mathfrak{g}$ with highest weight $\lambda \in P^{+}$, 
which can be regarded as a $\hat{\mathfrak{g}}_{\geq 0}$-module by letting $\hat{\mathfrak{g}}_{+}$ act trivially and the central element $\underline c$ act by a scalar $k$. 
The induced module 
$$V_{\lambda,k}=\text{Ind}_{\hat{\mathfrak{g}}_{\geq 0}}^{\hat{\mathfrak{g}}}V_{\lambda} \cong U(\hat{\mathfrak{g}}_{-}) \otimes V_{\lambda}$$ 
is called the Weyl module. 
It is $\mathbb{N}$-graded with $V_{\lambda, k} [0] = V_{\lambda}$. 

The following is well known (cf. [FZ]).

\begin{prop}\label{VOAaffine}
Let $k\neq -h^{\vee}$. 
The vacuum module $V_{0,k} \cong U(\hat{\mathfrak{g}}_{-}) \mathbf{1}$ 
is a vertex operator algebra of rank $\frac{k \cdot \text{dim} \mathfrak{g}}{k+h^{\vee}}$ 
with the Virasoro element given by $ \omega = \frac{1}{2(k+h^{\vee})} $ $\sum u_{i}(-1)u_{i}(-1) \mathbf{1}$, 
where $ \{u_{i}\} $ is an orthonormal basis of $\mathfrak{g}$ with respect to  $(,)$. 
The Weyl modules $V_{\lambda, k} $'s become representations of $V_{0, k}$ and in particular
$L_{0}|_{V_{\lambda, k} [n]} = n+ \triangle (\lambda)$ where $\triangle (\lambda) = \frac{ ( \lambda, \lambda +2 \rho )}{2 (k+h^{\vee})}$.  
\end{prop}

The vertex operator algebra $V_{0, k}$ is generated by the quantum fields:
$$Y( x(-1) \mathbf{1}, z) = x(z) = \sum_{n \in \mathbb{Z}} x(n) z^{-n-1} $$ 
 with the operator product expansions (OPE): 
$$x(z_{1}) y(z_{2}) \sim \frac{ [x, y] (z_{2})}{ z_{1} -z_{2}} + \frac{k (x, y)} {(z_{1}-z_{2})^{2}}$$
for $x, y\in \mf{g}$. 

Recall that the \emph{normally ordered product} of two fields 
$Y(u, z_{1}) = \sum u_{(n)} z_{1}^{-n-1}$ and 
$Y(v, z_{2})=\sum v_{(n)} z_{2}^{-n-1}$ is defined to be
$$:Y(u, z_{1})Y(v, z_{2}): = Y_{-}(u, z_{1})Y(v, z_{2}) + Y(v, z_{2})Y_{+}(u, z_{1})$$
where
$Y_{-}(u, z_{1})=\sum _{n <0} u_{(n)} z_{1}^{-n-1}$, $Y_{+}(u, z_{1})= \sum _{n \geq 0} u_{(n)} z_{1}^{-n-1} $.
Then
$$Y(\omega, z) = \frac{1}{2 (k+h^{\vee})} \sum _{i} :u_{i}(z) u_{i}(z):$$
$$L_{n}= \omega_{(n+1)}= \frac{1}{2 (k+h^{\vee})} \sum _{i} \sum _{j \in \mathbb{Z}} :u_{i}(-j) u_{i} (j+n): $$ where $:u_{i}(p)u_{i}(q):$ stands for $u_{i}(p)u_{i}(q)$ if $p <q$ or $u_{i}(q)u_{i}(p)$ if $ p \geq q$.

\begin{lemma}
For $k \notin \mathbb{Q}$, the Weyl module $V_{\lambda, k}$ is irreducible for any $\lambda\in P^{+}$. 
\end{lemma}

\begin{proof}
See Corollary 2.8.3 in [EFK] for example.
\end{proof}

\begin{prop} \label{intertwiners}
Let $k \notin \mathbb{Q}$ and  $ \lambda, \mu, \nu \in P^{+}$. 
For any $f \in \text{Hom} _{\mf{g}} ( V_{\lambda} \otimes V_{\mu}, V_{\nu})$, 
there exists a unique intertwining operator
$$\Phi^{f}(\cdot, z): V_{\lambda, k} \to  \text{End}\, ( V_{\mu, k} , V_{\nu, k} ) \{ z \}$$ 
such that it can be written as
$$\Phi^{f}(v, z) = \sum _{ n\in \mathbb{Z}} v_{(n)} z^{-n-1} z^{ \triangle (\nu) - \triangle (\lambda)- \triangle (\mu)}$$
where $v_{(n)} V_{\mu, k} [m] \subset V_{\nu, k} [-n-1+\text{deg}\, v+m]$
for homogeneous $v\in V_{\lambda, k}$ and 
$ {v_{1}}_{(-1)} u_{1} = f (v_{1} \otimes u_{1})$ 
for any $v_{1}\in V_{\lambda} $, $u_{1}\in V_{\mu}$.
Hence
$\mathcal{V}_{V_{\lambda, k} V_{\mu, k}}^{V_{\nu, k}} \cong \text{Hom} _{\mf{g}} ( V_{\lambda} \otimes V_{\mu}, V_{\nu})$.
\end{prop}

\begin{proof}
See Theorem 3.1.1 in [EFK] or Theorem 1.5.3  in [FZ]  for example.
\end{proof}

\begin{corollary}\label{propertyofintertwiners}
Let $k, \lambda, \mu, \nu, f, \Phi$ be as in Proposition 2.6, then for any $v \in V_{\lambda, k}, v_{1}\in V_{\lambda}$, $x\in \mf{g}$, $m\in \mathbb{Z}$, $n\in \mathbb{N}$ we have 
$$[ x(m), \Phi^{f}(v, z)] = \sum _{l \geq 0} {m \choose l} \Phi^{f}( x(l) v, z) z^{m-l} $$ 
$$\Phi^{f}( x(-n-1) v, z) = :\partial _{z}^{(n)} x(z)\Phi^{f}(v, z):$$
where $\partial_{z}^{(n)} = \frac{1}{n!} \frac{d^{n}} {d z^{n}}$ and in particular
$$[ x(m), \Phi^{f}(v_{1}, z)] = \Phi^{f} (x \cdot v_{1}, z) z^{m}$$
$$[L_{-1},  \Phi^{f}(v, z)] =  \Phi^{f}(L_{-1}v, z) = \partial_z \Phi^{f} (v, z). $$
Let $\varkappa = k+h^{\vee}$, then  
$$\frac{d}{dz} \Phi^{f} (v_{1}, z)=  \frac{1}{\varkappa} \sum_{i} :u_{i}(z) \Phi^{f} (u_{i} \cdot v_{1}, z): .$$
\end{corollary}

\begin{proof}
Everything follows from the Jacobi identity except the last equation, for which see Theorem 3.2.2  in [EFK].
\end{proof}

\begin{notation}
For $\lambda, \mu, \nu, \gamma  \in P^{+}$, 
we denote $\text{Hom}_{\mf{g}} ( V_{\lambda} \otimes V_{\mu}, V_{\nu})$ by $V_{\lambda \mu}^{\nu}$. Similarly set $V_{\lambda ^{\ast} \mu ^{*}} ^{\nu ^{*}} = \text{Hom} _{\mf{g}} (V_{\lambda ^{*}} \otimes V_{\mu ^{*}} , V_{\nu ^{*}})$, 
$ V_{\lambda \mu \nu}^{\gamma} = \text{Hom} _{\mf{g}} (V_{\lambda} \otimes V_{\mu} \otimes V_{\nu} , V_{\gamma}) $ and etc. 
If $V$ is a $\mf{g}$-module, define $V^{\mf{g}} = \{ v\in V: x \cdot v=0   \,\,\,\,\,  \forall x \in \mf{g} \}$. 
\end{notation}

\begin{lemma}
One has 
$V_{\lambda \mu}^{\nu} \cong {(V_{\lambda}^{*} \otimes V_{\mu}^{*} \otimes V_{\nu} ) }^{\mf{g}}$, 
 $V_{\lambda ^{\ast} \mu ^{*}} ^{\nu ^{*}} \cong {(V_{\lambda} \otimes V_{\mu} \otimes V_{\nu}^{*})}^{\mf{g}} $,
 $ V_{\lambda \mu \nu}^{\gamma} \cong { (V_{\lambda}^{*} \otimes V_{\mu}^{*} \otimes V_{\nu}^{*} \otimes V_{\gamma} )}^{\mf{g}}$ 
 and etc. 
\end{lemma}

\begin{proof}
The isomorphisms are canonical.
\end{proof}

Clearly $f \mapsto f^{*}$ establishes an isomorphism between $\text{Hom} _{\mf{g}} ( V_{\lambda} \otimes V_{\mu}, V_{\nu})$  and $ \text{Hom}_{\mf{g}} ( V_{\nu}^{*}, V_{\lambda}^{*} \otimes V_{\mu}^{*})$. Let $g\in V_{\lambda ^{\ast} \mu ^{*}} ^{\nu ^{*}}$, then $g^{*} \in \text{Hom}_{\mf{g}}( V_{\nu}, V_{\lambda} \otimes V_{\mu})$. We define a bilinear form
$$(,): V_{\lambda \mu}^{\nu} \times V_{\lambda ^{\ast} \mu ^{*}} ^{\nu ^{*}} \to \mathbb{C} $$ by
$(f, g) = \frac{ \textrm{tr} (f g^{*}) }{ \textrm{dim} V_{\nu} }$. In fact $f g^{*} = (f, g) \textrm{Id}$ on $V_{\nu}$ since $\textrm{Hom}_{\mf{g}}( V_{\nu}, V_{\nu}) = \mathbb{C}$. Clearly $(,)$ is non-degenerate and by the identification made in previous lemma, $(,)$ is the same as the tensor product of natural pairings between $V_{\lambda}$ and $V_{\lambda}^{*}$, $V_{\mu}$ and $V_{\mu}^{*}$, $V_{\nu}$ and $V_{\nu}^{*}$. 

\begin{lemma}\label{bilinearform}
One has 
$\bigoplus _{\alpha} V_{\lambda \alpha}^{\gamma} \otimes V_{\mu \nu}^{\alpha} \cong V_{\lambda \mu \nu}^{\gamma}$,
$\bigoplus _{\beta} V_{\lambda ^{*} \beta ^{*}}^{\gamma ^{*}} \otimes V_{\mu ^{*} \nu ^{*}}^{\beta ^{*}} \cong V_{\lambda ^{*} \mu^{*} \nu ^{*}} ^{\gamma ^{*}}$. 
Under these isomorphisms, the pairing between 
$\bigoplus _{\alpha} V_{\lambda \alpha}^{\gamma} \otimes V_{\mu \nu}^{\alpha}$ and 
$\bigoplus _{\beta} V_{\lambda ^{*} \beta ^{*}}^{\gamma ^{*}} \otimes V_{\mu ^{*} \nu ^{*}}^{\beta ^{*}} $
($0$ if $\alpha \neq \beta$) 
is equivalent to the pairing between  
$V_{\lambda \mu \nu}^{\gamma}$ and $ V_{\lambda ^{*} \mu^{*} \nu ^{*}} ^{\gamma ^{*}}$. 
\end{lemma}

\begin{proof}
This is an easy exercise as well. 
\end{proof}

Let $\mathcal{R}(G)$ be the space of regular functions on the algebraic group $G$. 
By the Peter-Weyl decomposition 
$\mathcal{R}(G) = \bigoplus _{\lambda \in P^{+}} V_{\lambda} \otimes V_{\lambda} ^{*}$
consists of the matrix coefficients of all irreducible finite dimensional representations of $G$. 
The following proposition explains how to express the (point-wise) multiplication of regular functions in terms of intertwining operators with respect to the left and right $\mf{g}$-actions. 

\begin{prop}\label{peterweyl}
Write $\mathcal{R}(G)= \bigoplus _{\lambda \in P^{+}} V_{\lambda} \otimes V_{\lambda} ^{*}$, 
then for any $u_{1} \otimes u_{1}' \in V_{\lambda} \otimes V_{\lambda} ^{*}$ 
and $v_{1} \otimes v_{1}' \in V_{\mu} \otimes V_{\mu} ^{*}$,  
$$(u_{1} \otimes  u_{1}') \cdot (v_{1} \otimes v_{1}') = \sum _{\nu \in P^{+}}  \sum _{i=1}^{dim V_{\lambda \mu}^{\nu}} f_{i}(u_{1} \otimes v_{1}) \otimes g_{i} ( u_{1}' \otimes v_{1}' ) $$
 where $\{f_i\}$ is a basis of  $V_{\lambda \mu}^{\nu}$ 
 and $\{g_i\}$ is the dual basis in $V_{\lambda ^{\ast} \mu ^{*}} ^{\nu ^{*}}$ 
 with respect to the bilinear form $(,)$ 
 between $V_{\lambda \mu}^{\nu}$ and  $V_{\lambda ^{\ast} \mu ^{*}} ^{\nu ^{*}}$. 
 \end{prop}

\begin{proof}
First it is clear that the right hand side does not depend on the choice of basis. 
Fix a decomposition $V_{\lambda} \otimes V_{\mu} = \oplus_{\nu} V_{\nu}^{ \oplus m_{\nu}}$ 
where $m_{\nu} = \text{dim}_{\mathbb{C}} V_{\lambda \mu}^\nu$. 
Let $\{p_{i}, i=1,2, \ldots, m_{\nu}\}$ be the projections of 
$V_{\lambda} \otimes V_{\mu}$ onto various copies of $V_{\nu}$ in it, 
and let $\{q_{j}, j=1,2, \ldots, m_{\nu}\}$ be the dual injections of 
$V_{\nu}$ into $V_{\lambda} \otimes V_{\mu}$ 
such that $ p_{i} q_{j} = \delta_{ij}$ on $V_{\nu}$ and 
$\sum _{i} q_{i}p_{i}$ equals the projection of 
$V_{\lambda} \otimes V_{\mu}$ onto its $\nu$ components. 
The $\nu$-component in the Peter-Weyl decomposition of  
$(u_{1} \otimes u_{1}') \cdot (v_{1} \otimes v_{1}')$ 
is  $\sum _{i} p_{i} (u_{1} \otimes v_{1}) \otimes q_{i}^{*} (u_{1}' \otimes v_{1}' )$ 
where $q_{i}^{*}\in V_{\lambda ^{*} \mu ^{*}}^{\nu ^{*}}$.  
Note that $\{p_{i}\}$ is a basis of $V_{\lambda \mu}^{\nu}$ 
and $\{q_{j}\}$ is a basis of $\text{Hom} _{\mf{g}} (V_{\nu}, V_{\lambda} \otimes V_{\mu})$. 
Since $p_{i} (q_{j}^{*})^{*} = p_{i} q_{j} = \delta_{ij}$, 
it follows that $(p_{i}, q_{j}^{*}) = \delta _{ij}$ 
which means that $\{q_{i}^{*}\}$ is the dual basis of $\{p_{i}\}$.
The proposition follows. 
\end{proof}

\begin{theorem} \label{mainthm1}
Let $\mf{g}$ be a simple complex Lie algebra and $\varkappa \notin \mathbb{Q}$. 
Set $k=\varkappa - h^{\vee}, \bar{k} = - \varkappa - h^{\vee}$, 
then 
$$\V= \bigoplus_{\lambda \in P^{+}}V_{\lambda,k} \otimes V_{\lambda^{*},\bar{k}}$$ 
is a vertex operator algebra of rank $2 \,\textrm{dim}(\mf{g})$. 
\end{theorem}

Here $V_{\lambda ^{*}, \bar{k}}$ means the induced representation of 
$V_{\lambda}^{*}$ of level $\bar{k}$. 
We treat $\lambda ^{*}$ as $- \omega _{0} \lambda$. 
In fact $V_{\lambda ^{*}, \bar{k}}$ is isomorphic to the contragredient module 
$V_{\lambda, \bar{k}}^{c}$ (see Definition \ref{contragredientmodule}) 
since it is irreducible when $\bar k \notin \mathbb{Q}$. 
The rest of Section 2 is devoted to the proof of the theorem.

$\V$ is obviously $\mathbb{N}$-graded. 
Choose the vacuum vector $\mathbf{1} = 1\otimes 1 \in V_{0, k} \otimes V_{0, \bar{k}}$.
We construct the vertex operators by pairing the intertwining operators
between $V_{0, k}$-modules and $V_{0, \bar k}$-modules. 
Specifically for any $u \otimes u' \in V_{\lambda, k} \otimes V_{\lambda^{*}, \bar{k}}$,  
$v \otimes v'  \in V_{\mu,k} \otimes V_{\mu^{*}, \bar {k}}$, 
define 
$$Y(u \otimes u', z) v \otimes v' = \sum _{\nu \in P^{+}} \sum _{i=1}^{dim V_{\lambda \mu}^{\nu}} \Phi ^{f_{i}} (u, z) v \otimes \Psi ^{g_{i}} (u' , z) v' $$ 
where
$\{f_{i}\}$ and $\{g_{i}\}$ are dual basis in $V_{\lambda \mu}^{\nu}$ 
and $V_{\lambda^{*} \mu^{*}}^{\nu^{*}} $, and 
$$\Phi^{f_{i}} (\cdot, z) : V_{\lambda, k} \otimes V_{\mu, k} \to V_{\nu, k} \{ z \}$$ 
$$\Psi ^{g_{i}} (\cdot, z) : V_{\lambda^{*}, \bar{k}} \otimes V_{\mu ^{*}, \bar{k}} \to V_{\nu ^{*}, \bar{k}} \{ z \}$$ 
are the associated intertwining operators (see Proposition \ref{intertwiners}). 
Since $\omega_{0}  \rho = - \rho$, we have 
$$\bar{\triangle} (- \omega_{0}  \lambda) =\frac { (-\omega_{0} \lambda, -\omega_{0}  \lambda+2\rho)} {2(\bar{k} + h^{\vee})} = \frac{( \lambda, \lambda+2\rho )}{-2\varkappa} = - \triangle(\lambda), $$ 
hence the $z^{\triangle (\nu) -\triangle (\lambda) -\triangle (\mu)}$ term 
in $\Phi ^{f_{i}}( \cdot, z)$ cancels with 
the $z^{\bar{\triangle} (\nu^{*}) - \bar{\triangle}(\lambda^{*}) - \bar{\triangle} (\mu^{*})}$ term 
in $\Psi^{g_{i}} (\cdot, z)$, 
therefore only terms of integral powers of $z$ remain. 
It is clear that the vertex operators of homogeneous elements satisfy the degree condition. 

The representation $Y_{V_{\lambda, k}} (\cdot, z):  V_{0, k} \otimes V_{\lambda, k} \to V_{\lambda, k} [[z^{\pm 1}]]$ is the intertwining operator corresponding to $1\in V_{0 \lambda}^{\lambda} \cong \mathbb{C}$, i.e. $Y_{V_{\lambda, k}}(\cdot, z) = \Phi^{1} (\cdot, z)$; similarly $Y_{V_{\lambda ^{*}, \bar{k}}}(\cdot, z) = \Psi ^{1}(\cdot, z)$. Hence 
$Y(\mathbf{1},z)u\otimes u'=Y_{V_{\lambda,k}}(1,z)u\otimes  Y_{V_{\lambda^{*},\bar{k}}}(1,z)u' =u\otimes u'$,
i.e. $Y(\mathbf{1},z) =\textrm{Id}$ and for any $x,y\in \mf{g}$, we have 
$Y(x(-1)\mathbf{1},z)u\otimes u'=x(z)u\otimes u'$ and 
$Y(\bar{y}(-1)\mathbf{1},z)u\otimes u'= u\otimes \bar{y}(z)u'$ 
where the 'bar' is used to denote the $\hat{\mf{g}}$-action of level $\bar{k}$. 
Indeed the vertex operators associated to elements in $V_{0, k}\otimes V_{0, \bar{k}}$ 
 generate the left and right $\hat{\mf{g}}$-actions on $\V$. 

Set 
$$\omega = \frac{1}{2\varkappa} \sum_{i} u_{i}(-1)u_{i}(-1) \mathbf{1} - \frac{1}{2\varkappa} \sum_{i} \bar{ u_{i}}(-1)\bar{u_{i}}(-1) \mathbf{1} $$ 
where $\{u_i\}$ is an orthonormal basis of $\mf{g}$ with respect to the normalized Killing form $(, )$ 
(see Proposition \ref{VOAaffine}) 
and let $Y(\omega, z) = \sum_{n\in \mathbb{Z}} \mathcal{L}_{n} z^{-n-2}$,
then $\mathcal{L}_{n} = L_{n} + \bar{L} _{n}$ where 
$$L_{n} =  \frac{1}{2\varkappa} \sum _{i} \sum _{j \in \mathbb{Z}} :u_{i}(-j) u_{i} (j+n): 
\quad \bar{L}_{n} =  \frac{1}{2(-\varkappa)} \sum_{i} \sum_{j \in \mathbb{Z}} :\bar{u_{i}}(-j) \bar{u_{i}} (j+n): .$$
Since $L_{n}, \bar{L}_{n}$ satisfy the Virasoro commutation relations of central charges 
$ \frac{k\cdot \text{dim} \mf{g}} {\varkappa}$ and $ \frac{\bar{k} \cdot \text{dim} \mf{g}} {-\varkappa}$ respectively and they commute with each other, 
$\mathcal{L}_{n}$ satisfy the Virasoro relations of central charge 
$$\frac{k\cdot \text{dim} \mf{g}}{\varkappa} +  \frac{\bar{k} \cdot \text{dim} \mf{g}} {-\varkappa} = \frac{ [(\varkappa - h^{\vee})- (-\varkappa - h^{ \vee})] \cdot \text{dim} \mf{g}} {\varkappa} = \frac{ 2 \varkappa \cdot \text{dim}\mf{g}} {\varkappa} = 2 \text{ dim} \mf{g}. $$ 

By Proposition \ref{VOAaffine}, 
we have 
$L_{0}| _{V_{\lambda, k}[n]} = n+\triangle (\lambda)$ and 
$\bar{L}_{0} |_{V_{\lambda^*, \bar k} [m]} = m+\bar{\triangle} ( \lambda^{*}) = m- \triangle (\lambda)$, 
hence $\mathcal{L} _{0} $ does act as the gradation operator on $\V$. 

It's known that 
$\mathcal{V} _{V_{0, k} V_{\lambda, k}} ^{V_{\lambda, k}} \cong \mathcal{V} _{V_{\lambda, k} V_{0, k}} ^{V_{\lambda, k}} $ (cf. [FHL]). 
The image of $Y_{V_{\lambda, k}} (\cdot, z)$ under the isomorphism is  $\Phi^{1}(\cdot, z)$ 
where $1\in V_{\lambda 0}^{\lambda} \cong \mathbb{C}$. 
More specifically 
$\Phi^{1} (u, z) v = e^{ z L_{-1}} Y_{V_{\lambda, k}} (v, -z) u$ 
for any $ v \in V_{0,k}, u\in V_{\lambda, k}$. 
Set $v=1$, one gets $\Phi^{1} (u, z) 1 =   e^{ z L_{-1}} u$; similarly 
$\Psi^{1} (u', z) 1= e^{z \bar{L}_{-1}} u'$ for any $u' \in V_{\lambda^{*}, \bar{k}}$. Hence 
$Y( u \otimes u', z) \mathbf{1} = \Phi^{1} ( u, z) 1 \otimes \Psi^{1} (u', z) 1= e^ {z L_{-1}} u \otimes e^{z \bar{L}_{-1}} u' = e^{z \mathcal{L} _{-1}} u \otimes u' $, 
in particular $Y(u \otimes u', z) \mathbf{1} |_{z=0} = u \otimes u'$. 

Define $\partial = \mathcal{L}_{-1}$. Since 
$[L_{-1}, \Phi^{f_{i}} (u, z)] = \Phi ^{f_{i}} ( L_{-1} u, z) = \frac{d}{dz} \Phi ^{f_{i}}(u, z)$ and 
$[\bar{L}_{-1}, \Psi^{g_{i}} (u', z)] = \Psi ^{g_{i}} (\bar{L}_{-1} u',  z) = \frac{d}{dz} \Psi ^{g_{i}}(u', z)$,
it follows  that
$$[\mathcal{L}_{-1}, Y(u \otimes u', z)] =Y (\mathcal{L}_{-1} (u \otimes u'), z) = \frac{d}{dz}Y(u \otimes u', z). $$

Now it remains to show that any two fields are mutually local. 
First note that by Corollary \ref{propertyofintertwiners}, 
for any $u_{1} \otimes u_{1}' \in V_{\lambda} \otimes V_{\lambda}^{*}$, $x, y\in \mf{g}$, 
we have
$$[x(m), Y(u_{1}\otimes u_{1}', z)] = z^{m} Y( x\cdot u_{1} \otimes u_{1}', z)$$
$$[\bar{y}(n), Y(u_{1}\otimes u_{1}', z)] = z^{n} Y( u_{1} \otimes \bar{y}\cdot u_{1}', z). $$
It implies that 
$(z_{1}-z_{2})[x(z_{1}), Y(u_{1}\otimes u_{1}', z_{2})]=(z_{1}-z_{2})[\bar{y}(z_{1}), Y(u_{1}\otimes u_{1}', z_{2})]=0$, i.e. $x(z)$, $\bar{y}(z)$ are local with $Y(a, z), a\in A$. 

Since $Y( x_{1}(-p_{1})\ldots x_{m}(-p_{m})u_{1} \otimes \bar{y}_{1}(-q_{1})\ldots \bar{y}_{n}(-q_{n}) u_{1}', z) $ can be expressed as normally ordered products of $x_{i}(z), \bar{y}_{j}(z)$ and
$Y(u_{1} \otimes u_{1}', z)$,  it suffices to prove the locality of fields associated to vectors in $\V_0$  due to the following lemma.

\begin{lemma}
If $A(z), B(z), C(z)$ are three mutually local fields, then the fields $:A(z)B(z):$ and $C(z)$ are also mutually local.
\end{lemma}

\begin{proof}
This is Proposition 5.5.15. in [LL] or Lemma 2.3.4. in [FB].
\end{proof}

Let $a, b\in \V_0$. If $\V$ is a vertex algebra, by consideration of grading, it must be true that 
$a_{(i)} b = 0 $ for any $i \geq 0$, which implies that  $[Y(a, z_{1}),Y(b, z_{2})]=0$ by Borcherds' commutator formula.  So our goal is to prove that
$Y(a, z_{1})$ and $ Y(b, z_{2})$ commute with each other.  
Since 
$$Y(a, z_{1}) Y(b, z_{2}) x(-m) v$$
$$
= x(-m)Y(a, z_{1}) Y(b, z_{2}) v- z_{1}^{-m} Y(x\cdot a, z_{1}) Y(b, z_{2}) v - z_{2}^{-m} Y(a, z_{1})Y(x\cdot b, z_{2}) v$$
and the parallel equalities hold for the other copy of the $\hat{\mf{g}}$-action, 
by induction on the degree of $v$, 
it suffices to prove that $Y(a, z_{1})Y(b, z_{2})c= Y(b, z_{2})Y(a, z_{1})c$ for any $c \in \V_0$.

 \begin{define} \label{contragredientmodule}
Let $M=\oplus_{n} M_{n}$ be a $\mathbb{Z}$-graded $\hat{\mf{g}}$-representation of central charge $k$ with $\text{dim}M_{n} < \infty$,  define the contragredient module $M^{c}$ as follows: as a vector space $M^{c} =\oplus_{n} M_{n}^{*}$; the $\hat{\mf{g}}$-action is given by $(x(l) v^{*}, v) = -(v^{*}, x(-l) v)$ for $x\in \mf{g}$, $l\in \mathbb{Z}$, $v^{*}\in M^{c}$, $v\in M$ and the center $\underline c$ acts by $k\, \text{Id}$. 
$M^{c}$ is $\mathbb{Z}$-graded as well. 
\end{define}

Let $\V^{c}= \oplus_{\lambda \in P^{+}, n \geq 0} (V_{\lambda, k} \otimes V_{\lambda^{*}, \bar{k}})[n]^{*}$
be the restricted dual space of $\V$. 
For $k, \bar{k} \notin \mathbb{Q}$, 
we have 
$(V_{\lambda, k})^{c}\cong V_{\lambda^{*}, k}$ and 
$(V_{\lambda^{*}, \bar{k}})^{c} \cong V_{\lambda, \bar k}$,  
hence 
$\V^{c}$ is generated by 
$\V_0^*= \oplus_{\lambda} (V_\lambda \otimes V_\lambda^*) ^*$ 
as a $\hat{\mf{g}} \oplus \hat{\mf{g}}$-module. 
Since 
$$(x(-m) v^{*}, Y(a, z_{1})Y(b, z_{2})c) =-(v^{*}, x(m)Y(a, z_{1})Y(b,z_{2})c)$$
$$
=-z_{1}^{m}(v^{*}, Y(x\cdot a, z_{1})Y(b, z_{2})c) -z_{2}^{m} (v^{*}, Y(a, z_{1})Y(x\cdot b, z_{2})c)
$$        
for any $x\in \mf{g}$, $m>0$, $v^{*}\in \V^c$ 
and the parallel equalities involving the other copy of the $\hat{\mf{g}}$-action hold, 
again by induction on the degree of $v^{*}$, 
it suffices to prove that
$(d^{*}, Y(a,z_{1})Y(b, z_{2})c)=(d^{*}, Y(b, z_{2}) Y(a, z_{1})c)$ 
for any $a, b, c\in \V_0$, $d^{*}\in \V_0^*$.
 
Suppose $a=u_{1} \otimes u_{1}' \in V_{\lambda}\otimes V_{\lambda}^{*}$, 
$b=v_{1} \otimes v_{1}' \in V_{\mu} \otimes V_\mu^*$, 
$c= w_{1} \otimes w_{1}' \in V_ \nu \otimes V_{\nu}^{*}$ and 
$d^{*} = d_{1}' \otimes d_{1} \in V_{\gamma}^{*} \otimes V_{\gamma}$,
then
$$(d^{*}, Y(a, z_{1})Y(b, z_{2})c)$$
$$=\sum_{ i, j,\alpha} (d_{1}' ,  \Phi^{r_{i}^{\alpha}}( u_{1}, z_{1} ) \Phi ^{f_{j}^{\alpha}}( v_{1}, z_{2}) w_{1}   ) (d_{1}, \Psi^{s_{i}^{\alpha}}( u_{1}', z_{1}) \Psi^{g_{j}^{\alpha}}(v_{1}', z_{2}) w_{1}')$$ 
where
$\{f_{j}^{\alpha} \}$ and $\{g_{j}^{\alpha}\}$ are dual basis in 
$V_{\mu \nu}^{\alpha}$ and $V_{\mu^{*} \nu^{*}}^{\alpha^{*}}$ 
while 
$\{r_{i}^{\alpha} \}$ and $\{s_{i}^{\alpha} \}$ are dual basis in 
$ V_{\lambda \alpha}^{\gamma} $ and  $ V_{\lambda^{*} \alpha^{*}}^{\gamma^{*}} $
for various $\alpha \in P^+$.  
By Lemma \ref{bilinearform}, 
$r_{i}^{\alpha} \otimes f_{j}^{\alpha} \in \oplus_{\alpha} V_{\lambda \alpha}^{\gamma} \otimes V_{\mu \nu}^{\alpha} \cong V_{\lambda \mu \nu}^{\gamma} $ and 
$s_{i}^{\alpha} \otimes g_{j}^{\alpha} \in \oplus_{\alpha} V_{\lambda^{*} \alpha^{*}}^{\gamma^{*}} \otimes V_{\mu^{*} \nu^{*}}^{\alpha^{*}} \cong V_{\lambda^{*} \mu^{*} \nu^{*}}^{\gamma^{*}} $ 
are dual basis as well,  i.e. 
$(r_{i}^{\alpha} \otimes f_{j}^{\alpha}, s_{i'}^{\beta} \otimes g_{j'}^{\beta}) = \delta _{\alpha \beta} \delta_{i i'} \delta_{j j'}$ and 
$$ \sum_{i,j,\alpha} (r_{i}^{\alpha} \otimes f_{j}^{\alpha} ) (u_{1}\otimes v_{1} \otimes w_{1}) \otimes (s_{i}^{\alpha}\otimes g_{j}^{\alpha} ) (u_{1}' \otimes v_{1}' \otimes w_{1}' ) \in V_{\gamma} \otimes V_{\gamma}^{*} $$ 
is the $\gamma$-component of $a\cdot b\cdot c$. 

Denote
$$\Phi^{i j \alpha} = (\cdot, \Phi^{r_{i}^{\alpha}}(\cdot, z_{1}) \Phi^{f_{j}^{\alpha}}(\cdot, z_{2}) \cdot )\in V_{\lambda \mu \nu}^{\gamma} \{z_{1}, z_{2}\}$$
$$\Psi^{i' j' \beta} = (\cdot, \Psi^{s_{i'}^{\beta}}(\cdot, z_{1}) \Psi^{g_{j'}^{\beta}} (\cdot, z_{2}) \cdot) \in V_{\lambda^{*} \mu^{*}\nu^{*}}^{\gamma^{*}} \{z_{1}, z_{2}\}, $$
then
$$ (\cdot _{\gamma}, Y(\cdot_{\lambda}, z_{1}) Y(\cdot_{\mu}, z_{2})  \cdot_{\nu}) = \sum_{i,j,\alpha} \Phi^{ij\alpha} \otimes \Psi^{ij\alpha}. $$

 
\begin{lemma} \label{215}
$\Phi^{i j \alpha} $ and $\Psi^{i' j' \beta}$ satisfy the following Knizhnik-Zamolodchikov equations:
$$\frac{d}{d{z_{1}}} \Phi^{ij\alpha} = \frac{1}{\varkappa}\Phi^{ij\alpha} (\frac{\Omega_{12}}{z_{1}-z_{2}} + \frac{\Omega_{13}}{z_{1}})$$
 $$\frac{d}{d z_{2}}\Phi^{ij\alpha} = \frac{1}{\varkappa}\Phi^{ij\alpha} (\frac{\Omega_{12}}{-z_{1}+z_{2}} +\frac{\Omega_{23}}{z_{2}})$$
and
 $$\frac{d}{d z_{1}}\Psi^{i' j' \beta} = -\frac{1}{\varkappa} \Psi^{i' j' \beta} (\frac{\Omega_{12}}{z_{1}-z_{2}} + \frac{\Omega_{13}}{z_{1}})$$
$$\frac{d}{d z_{2}}\Psi^{i' j' \beta} = -\frac{1}{\varkappa} \Psi^{i' j' \beta} (\frac{\Omega_{12}}{-z_{1}+z_{2}} +\frac{\Omega_{23}}{z_{2}})$$
where $1/(z_{1}-z_{2})$ and $1/(-z_{1}+z_{2})$ stand for the formal power series expansion in $z_{2}/z_{1}$, and $\Omega_{12}$ acts on $V_{\lambda}\otimes V_{\mu} \otimes V_{\nu}$ by 
 $\sum_{i} u_{i}\otimes u_{i} \otimes 1$ (here $\{u_i\}$ is an orthonormal basis of $\mf{g}$ with respect to the normalized Killing form as before) and  $\Omega_{13}, \Omega_{23}$ act accordingly. 
 \end{lemma}  
 
 \begin{proof}
 See Theorem 3.4.1 in [EFK].
 \end{proof}
 
 The Knizhnik-Zamolodchikov equations can also be interpreted as identities of operator-valued analytic functions. In particular on the simply connected region $\mathcal{D} = \{ z_{1}, z_{2}\in \mathbb{C}: |z_{1}|>|z_{2}|>0, z_{1}, z_{2}\notin \mathbb{R}_{-} \}$, $ \Phi^{ij\alpha}$ (resp. $\Psi^{i' j' \beta}$) converge and span the space of solutions to the KZ equations on $\mathcal{D}$. 

 \begin{lemma} \label{216}
 For any $f\in V_{\lambda \mu \nu}^{\gamma}$ and  
$ g\in V_{\lambda^{*}  \mu^{*} \nu^{*}}^{\gamma^{*}}$, one has 
$ (f \Omega_{12}, g) = (f, g \Omega_{12})$, $(f \Omega_{13}, g) = (f, g \Omega_{13})$
and $(f \Omega_{23}, g) = (f, g \Omega_{23})$. 
 \end{lemma}
 
\begin{proof}
 It suffices to prove one of the equalities. It can be easily checked that 
 $\Omega_{12}^{*} = \Omega_{12}$, hence by definition $(f, g\Omega_{12}) = \frac{\text{tr} (f (g\Omega_{12})^{*})}{\text{dim}V_{\gamma}} =
\frac{\text{tr}( f\Omega_{12} g^{*})}{\text{dim}V_{\gamma}} = (f\Omega_{12}, g)$. 
\end{proof}
 
\begin{prop}
$(\Phi^{i j \alpha} , \Psi^{i' j' \beta}) = \delta_{ii'}\delta_{jj'}\delta_{\alpha \beta}$. 
\end{prop}
 
\begin{proof}
Since 
$\Phi^{i j \alpha} = z_{1}^{ \triangle(\gamma) -\triangle(\lambda)-\triangle(\alpha)} 
z_{2}^{ \triangle(\alpha)-\triangle(\mu)-\triangle(\nu)} ( r_{i}^{\alpha} \otimes f_{j}^{\alpha} + $ higher degree  terms in  $\frac{z_{2}}{z_{1}}) $ and 
$\Psi^{i' j' \beta} = z_{1}^{- \triangle(\gamma) +\triangle(\lambda)+\triangle(\beta)} z_{2}^{ -\triangle(\beta)+\triangle(\mu)+\triangle(\nu)} (s_{i'}^{\beta} \otimes g_{j'}^{\beta} +$ higher degree terms in $\frac{z_{2}}{z_{1}}) $, 
$(\Phi^{i j \alpha} , \Psi^{i' j' \beta}) \in  \mathbb{C} \{z_{1}, z_{2}\}$ is well defined. 
 By Lemma \ref{215} and Lemma \ref{216}, 
$$\frac{d}{d z_{1}}(\Phi^{i j \alpha} , \Psi^{i' j' \beta}) = (\frac{d}{d z_{1}} \Phi^{i j \alpha} , \Psi^{i' j' \beta}) + (\Phi^{i j \alpha} , \frac{d}{d z_{1}}\Psi^{i' j' \beta})$$
$$= \frac{1}{\varkappa}(\Phi^{ij\alpha} (\frac{\Omega_{12}}{z_{1}-z_{2}} + \frac{\Omega_{13}}{z_{1}}), \Psi^{i' j' \beta}) -\frac{1}{\varkappa}(\Phi^{i j \alpha} ,   \Psi^{i' j' \beta} (\frac{\Omega_{12}}{z_{1}-z_{2}} + \frac{\Omega_{13}}{z_{1}})) =0 $$
and similarly $\frac{d}{d z_{2}}(\Phi^{i j \alpha} , \Psi^{i' j' \beta})=0$. Hence $(\Phi^{i j \alpha} , \Psi^{i' j' \beta})$ is a constant, and this constant is nothing else but 
$ (r_{i}^{\alpha} \otimes f_{j}^{\alpha}, s_{i'}^{\beta} \otimes g_{j'}^{\beta} ) = \delta_{ii'}\delta_{jj'}\delta_{\alpha \beta}$.
\end{proof}

Think of $ \Phi^{ij\alpha} $ and $ \Psi^{i'j'\beta}$ as operator-valued ($ V_{\lambda \mu \nu}^{\gamma}$ and $V_{\lambda^{*} \mu^{*}\nu^{*}}^{\gamma^{*}}$-valued) functions on $\mathcal{D}$. 
At any point $(z_{1}, z_{2})\in \mathcal{D}$,  
one has $(\Phi^{ij\alpha} (z_{1} , z_{2}), \Psi^{i'j'\beta}(z_{1}, z_{2})) $ $= \delta_{ii'} \delta_{jj'} \delta_{\alpha \beta}$, %
which implies that $\{ \Phi^{ ij\alpha}(z_{1}, z_{2})\}$ constitute a basis  of
$V_{\lambda \mu \nu}^{\gamma}$ and 
$\{ \Psi^{ij\alpha}(z_{1}, z_{2})\}$ is the dual basis in $V_{\lambda^{*} \mu^{*} \nu^{*}}^{\gamma^{*}}$. 
Therefore
$$\sum_{i,j,\alpha} \Phi^{i j \alpha} \otimes \Psi^{i j \alpha}|_{(z_{1}, z_{2})} \equiv \sum_{i, j, \alpha} (r_{i}^{\alpha}\otimes f_{j}^{\alpha}) \otimes (s_{i}^{\alpha} \otimes g_{j}^{\alpha}), $$
 the canonical element in 
$V_{\lambda \mu \nu}^{\gamma} \otimes V_{\lambda^{*} \mu^{*}\nu^{*}}^{\gamma^{*}}$ 
associated with the non-degenerate bilinear form $(, )$ between 
$V_{\lambda \mu \nu}^{\gamma}$ and $V_{\lambda^{*} \mu^{*}\nu^{*}}^{\gamma^{*}}$. 
Hence 
$$(d^{*}, Y(a, z_{1}) Y(b, z_{2}) c) = (d^{*}, \sum_{i,j,\alpha} \Phi^{ij\alpha} (u_{1}\otimes v_{1} \otimes w_{1})\otimes \Psi^{ij\alpha} (u_{1}' \otimes v_{1}' \otimes w_{1}'))$$
$$\equiv  (d^{*}, \sum_{i,j,\alpha} (r_{i}^{\alpha} \otimes f_{j}^{\alpha} ) (u_{1}\otimes v_{1} \otimes w_{1}) \otimes (s_{i}^{\alpha}\otimes g_{j}^{\alpha} ) (u_{1}' \otimes v_{1}' \otimes w_{1}' )) = (d^{*}, a\cdot b\cdot c),$$
hence 
$(d^{*}, Y(a, z_{1})Y(b, z_{2})c )= (d^{*}, Y(b, z_{2}) Y(a, z_{1}) c)$ 
because $ a\cdot b \cdot c = b \cdot a \cdot c$. 
Theorem \ref{mainthm1} is now proved.

\section{Construction from vertex algebroid}
Let us recall how Gorbounov, Malikov and Schechtman construct the 
enveloping algebra of a vertex algebroid in [GMS1]. We fix the ground field to be $\mathbb{C}$ and adopt their notations for the preliminary. 


\begin{define}
A \emph{vertex algebroid} is a septuple 
$\mathcal{A} = (A, T, \Omega, \partial, \gamma, \langle, \rangle, c)$ 
where $A$ is a commutative algebra, 
$T$ is a Lie algebra acting by derivations on $A$ and equipped with the structure of an $A$-module, 
$\Omega$ is an $A$-module equipped with the structure of a module over the Lie algebra $T$, 
$\partial: A \to \Omega$ is an $A$-derivation and a morphism of $T$-modules, 
$\langle, \rangle:(T\oplus \Omega) \times (T\oplus \Omega) \to A$ is a symmetric bilinear pairing which is zero on $\Omega \times \Omega$, 
$c: T\times T \to \Omega$ is a skew symmetric bilinear pairing and 
$\gamma: A\times T \to \Omega$ is a bilinear map 
such that the following properties are satisfied for any $a, b\in A$, 
$\tau, \tau_{i}, \nu \in T$ and 
$\omega \in \Omega$:
$$[\tau, a\nu] = a [\tau, \nu] + \tau(a) \nu$$
$$ (a\tau)(b)= a \tau(b)$$
$$\langle \tau, \partial a \rangle = \tau(a)$$
$$\tau (a \omega) = \tau(a) \omega + a \tau(w)$$
$$ (a\tau) (\omega) = a\tau(\omega) + \langle \tau, \omega \rangle  \partial a$$
$$\tau( \langle \nu, \omega \rangle) = \langle  [\tau, \nu], \omega \rangle + \langle \nu, \tau(\omega) \rangle$$
$$\gamma(a, b\tau) = \gamma (ab, \tau) - a \gamma(b, \tau) - \tau(a) \partial b- \tau(b)\partial a$$
$$\langle a\tau_{1}, \tau_{2} \rangle= a\langle \tau_{1}, \tau_{2} \rangle + \langle \gamma(a, \tau_{1}), \tau_{2} \rangle -\tau_{1}\tau_{2}(a)$$
$$ c(a\tau_{1}, \tau_{2}) = ac(\tau_{1}, \tau_{2}) + \gamma(a, [\tau_{1}, \tau_{2}]) \gamma(\tau_{2}(a), \tau_{1}) + \tau_{2}(\gamma(a, \tau_{1}))$$
$$-\frac{1}{2}\langle \tau_{1}, \tau_{2} \rangle \partial a +\frac{1}{2}\partial \tau_{1}\tau_{2}(a)-\frac{1}{2}\partial \langle \tau_{2}, \gamma(a, \tau_{1}) \rangle $$
$$\langle [\tau_{1}, \tau_{2}], \tau_{3} \rangle + \langle \tau_{2}, [\tau_{1}, \tau_{3}] \rangle = \tau_{1}(\langle \tau_{2}, \tau_{3}\rangle ) -\frac{1}{2}\tau_{2}(\langle \tau_{1}, \tau_{3} \rangle ) -\frac{1}{2} \tau_{3}(\langle \tau_{1}, \tau_{2} \rangle)$$
$$+ \langle \tau_{2}, c(\tau_{1}, \tau_{3} )\rangle + \langle \tau_{3}, c(\tau_{1}, \tau_{2} )\rangle$$
$$d_{\textrm{Lie}} c(\tau_{1},\tau_{2},\tau_{3}) =-\frac{1}{2} \partial \{\langle[\tau_{1}, \tau_{2}], \tau_{3}\rangle +\langle[\tau_{1}, \tau_{3}],\tau_{2}\rangle-\langle[\tau_{2}, \tau_{3}],\tau_{1} \rangle-\tau_{1}(\langle \tau_{2}, \tau_{3} \rangle)$$
$$+\tau_{2}(\langle \tau_{1}, \tau_{3}\rangle) -2 \langle\tau_{3}, c(\tau_{1}, \tau_{2})\rangle \}$$
where $d_{\textrm{Lie}} c(\tau_{1},\tau_{2},\tau_{3}) =
\tau_{1} c(\tau_{2}, \tau_{3}) - \tau_{2} c(\tau_{1}, \tau_{3}) + \tau_{3} c(\tau_{1}, \tau_{2}) - c([\tau_{1}, \tau_{2}], \tau_{3})+ c([\tau_{1}, \tau_{3}], \tau_{2}) - c([\tau_{2}, \tau_{3}], \tau_{1})$. 
(cf. [GMS1] 1.1, 1.4)
\end{define}

The following proposition tells that to any vertex algebra with a splitting one can canonically associate a vertex algebroid. In fact the above axioms come from the Jacobi identity of a vertex algebra. 


\begin{prop}\label{VOAtoVA}
Let $V$ be an $\mathbb{N}$-graded vertex algebra over $\mathbb{C}$. 
Set $A= V_{0}$, then $A$ becomes a commutative algebra with the multiplication induced by the $_{(-1)}$ operation on $A$ and the vacuum vector becomes the identity. 
Let $\Omega \subset V_{1}$ be the $\mathbb{C}$-subspace generated by elements $a_{(-1)}\partial b$, $a, b\in A$ 
and let $\partial: A \in \Omega$ be the restriction of $\partial$ to $A$. 
Set $T=V_{1}/ \Omega$
and let $\pi: V_{1} \to T$ be the projection. 
Suppose $s: T \to V_{1}$ is a linear map such that $\pi \circ s = \text{Id}_T$.
Define a skew symmetric operation $[,]:V_{1} \times V_{1} \to V_{1}$ 
by $[x,y] =\frac{1}{2} (x_{(0)}y-y_{(0)}x)$ and set
$$ \gamma(a, \tau) = s(a\tau)-as(\tau),$$
$$\langle \tau_{1}+\omega_{1}, \tau_{2}+\omega_{2} \rangle = (s(\tau_{1})+\omega_{1})_{(1)} (s(\tau_{2})+\omega_{2}),$$
$$c(\tau_{1}, \tau_{2}) = s([\tau_{1}, \tau_{2}]) -[s(\tau_{1}), s(\tau_{2})].$$
Then the septuple $\mathcal{A} = (A, T, \Omega, \partial, \gamma, \langle, \rangle, c)$ 
is a vertex algebroid.
\end{prop}

\begin{proof}
This is Theorem 2.3 in [GMS1].
\end{proof}

We are most interested in the following example of a vertex algebroid.


\begin{prop} \label{11}
Let $A = \mathcal{R} (G) $ be the space of regular functions on a simple complex Lie group $G$ with Lie algebra $\mf{g}$, 
$T \cong A\otimes \mf{g}$ be the Lie algebra of vector fields over $G$ where $\mf{g}$ act on $A$ as left invariant vector fields, 
$\Omega:= \text{Hom}_{A} (T, A)$ be 1 forms, 
$\partial: A\to \Omega$ be the canonical differential given by $\partial a(\tau) =\tau(a)$ for $\tau \in T$, 
$\langle, \rangle : T \times \Omega \to A$ be the canonical pairing, 
$(, )$ be the normalized invariant form on $\mf{g}$ and 
$k$ be any complex number.  
The zero map $A \times \mf{g} \to \Omega$
(resp. $k(,) : \mf{g} \times \mf{g} \to \mathbb{C}$ 
and the zero map $\mf{g} \times \mf{g} \to \Omega$) 
can be extended to $\gamma : A \times T \to \Omega$
(resp. $\langle, \rangle: T\times T \to A$ and $c: T\times T \to \Omega$) 
in accordance with the axioms in Definition 3.1. 
Then $\mathcal{A}_{\mf{g},k} = (A, T, \Omega, \partial, \gamma, \langle, \rangle, c)$ is a vertex algebroid. 
\end{prop}

\begin{proof}
This is Example 1.11 which follows from Theorem 1.9 (Extension of Identities) in [GMS1].
\end{proof}

\begin{define}
A 1-truncated vertex algebra is a septuple 
$v = (V_{0}, V_{1}, \mathbf{1}, \partial,$ $ _{(-1)}, _{(0)}, _{(1)})$ 
where $V_{0}, V_{1}$ are $\mathbb{C}$-vector spaces, 
$\mathbf{1} \in V_{0}$ (the vacuum), 
$\partial : V_{0} \to V_{1}$ is a linear map and 
$_{(i)}: (V_{0} \oplus V_{1}) \times (V_{0} \oplus V_{1}) \to (V_{0} \oplus V_{1}) (i=-1, 0, 1)$ 
are bilinear operations of degree $-i-1$ 
such that the following axioms are satisfied for any $a,b,c\in V_{0}, x, y \in V_{1}$:
$$a_{(-1)}\mathbf{1} = a \qquad x_{(-1)} \mathbf{1} =x \qquad x_{(0)} \mathbf{1} = 0$$
$$(\partial a)_{(0)}b=0 \qquad (\partial a)_{(0)} x=0 \qquad (\partial a)_{(-1)} x= - a_{(0)}x$$
$$\partial ( a_{(-1)} b)=(\partial a)_{(-1)} b+ a_{(-1)} \partial b \qquad \partial(x_{(0)}a)=x_{(0)} \partial a$$
$$a_{(-1)}b =b_{(-1)} a \qquad a_{(-1)}x = x_{(-1)} a -\partial (x_{(0)} a)$$
$$x_{(0)}a = -a_{(0)} x \qquad  x_{(0)}y = -y_{(0)} x +\partial (y_{(1)}x)$$
$$ x_{(1)} y =y_{(1)} x$$
$$(a_{(-1)}b)_{(-1)}c=a_{(-1)}b_{(-1)}c$$
$$\alpha_{(0)}\beta_{(i)}\gamma=(\alpha_{(0)}\beta)_{(i)}\gamma+\beta_{(i)}\alpha_{(0)}\gamma \qquad (\alpha,\beta,\gamma \in V_{0}\oplus V_{1})$$
whenever the both sides are defined and
$$(a_{(-1)}x)_{(0)}b=a_{(-1)}x_{(0)}b$$
$$(a_{(-1)}b)_{(-1)}x=a_{(-1)}b_{(-1)}x+(\partial a)_{(-1)}b_{(0)}x+(\partial b)_{(-1)} a_{(0)}x$$
$$(a_{(-1)}x)_{(1)}y=a_{(-1)}x_{(1)}y-x_{(0)}y_{(0)}a. $$
(cf. [GMS1] 3.1)
\end{define}

One can assign a 1-truncated vertex algebra $(V_{0}, V_{1},\ldots)$ to any $\mathbb{N}$-graded vertex algebra $V$ by simple truncation. 
In fact the above axioms for 1-truncated vertex algebras are derived from the skew symmetry and the Jacobi identity of the vertex algebras. 
To any (split) 1-truncated vertex algebra one can associate a vertex algebroid by the same construction as in Proposition \ref{VOAtoVA}. Conversely as the next proposition states, a vertex algebroid yields a 1-truncated vertex algebra. 


\begin{prop} \label{broidtotrun}
Let $\mathcal{A} =(A, T, \Omega, \ldots)$ be a vertex algebroid. 
Set $V_{0}=A$, 
$V_{1} = T\oplus \Omega$
and let $\partial: V_{0}\to V_{1}$ be the composition 
of $\partial :A\to \Omega$ with the obvious embedding $\Omega \hookrightarrow V_{1}$. 
Let $\mathbf{1} = 1\in A$ and define the operations $_{(-1)},_{(0)},_{(1)}$ as follows:
$$ a_{(-1)} b=ab \qquad a_{(-1)} \omega = a\omega \qquad a_{(-1)} \tau= a\tau-\gamma (a, \tau)$$
$$a_{(0)}b=a_{(0)}\omega =\omega_{(0)} \omega'=0 \qquad \tau_{(0)} a=\tau(a) \qquad \tau_{(0)}\omega =\tau(\omega)$$
$$\tau_{(0)}\tau' = [\tau,\tau'] - c(\tau,\tau') +\frac{1}{2} \partial \langle \tau, \tau' \rangle$$
$$x_{(1)} y =\langle x, y \rangle$$
for any $a, b\in A$, $\tau, \tau'\in T$, $\omega, \omega'\in \Omega$ and $x, y\in T\oplus \Omega$.
Then $(V_{0}, V_{1}, \mathbf{1}, \partial, _{(-1)}, \\ _{(0)}, _{(1)})$ is a 1-truncated vertex algebra. %
\end{prop}

\begin{proof}
See 3.3 in [GMS1].
\end{proof}

Forget the $_{(-1)}$ operation in a 1-truncated vertex algebra, one gets what is called a 1-truncated conformal algebra. 


\begin{define}
A 1-truncated conformal algebra is a quintuple $c = (C_{0}, C_{1}, \\ \partial,  _{(0)},  _{(1)})$ %
where $C_{(i)}$ are $\mathbb{C}$-vector spaces, 
$\partial: C_{(0)} \to C_{(1)}$ is a linear map 
and $_{(i)}: (C_{0} \oplus C_{1}) \times (C_{0}\oplus C_{1}) \to C_{0}\oplus C_{1}$ are bilinear operations of degree $-i-1$ 
such that the following axioms are satisfied:
$$(\partial a)_{(0)} =0 \qquad (\partial a)_{(1)} = - a_{(0)} \qquad \partial(x_{(0)} a)=x_{(0)} \partial a$$
$$ x_{(0)} a = -a_{(0)} x \qquad x_{(0)}y=-y_{(0)}x+\partial(y_{(1)} x) \qquad x_{(1)}y=y_{(1)}x$$
$$\alpha_{(0)}\beta_{(i)}\gamma =(\alpha_{(0)}\beta)_{(i)}\gamma+\beta_{(i)}\alpha_{(0)}\gamma$$
for any $a\in C_{0}, x, y\in C_{1}, \alpha, \beta,\gamma\in C_{0}\oplus C_{1}$. (cf. [GMS1] 9.6)
\end{define}


\begin{define}
An $\mathbb{N}$-graded conformal algebra is an $\mathbb{N}$-graded $\mathbb{C}$-vector space $C=\oplus_{i\geq 0} C_{i}$ 
with a linear map $\partial:C\to C$ of degree 1 and a family of bilinear operations 
$$_{(n)} : C \times C \to C;  \quad (a, b) \mapsto a_{(n)} b $$
of degree $-n-1$ for $n\geq 0$ such that
$$(\partial a)_{(n)} b = -n a_{(n-1)} b$$
$$a_{(n)} b = (-1)^{n+1} \sum_{j=0}^{\infty} (-1)^{j} \partial^{(j)} (b_{(n+j)} a)$$
$$a_{(m)}b_{(n)}c=b_{(n)}a_{(m)}c+\sum_{j=0}^{m} {m \choose j} (a_{(j)}b)_{(m+n-j)}c$$
where $\partial^{(j)} = \partial^{j} / j!$ for all $a,b,c\in C, m, n\in \mathbb{N}$. 
(cf. [GMS1] 0.4)
\end{define}

The following proposition describes how to extend a 1-truncated conformal algebra to an $\mathbb{N}$-graded conformal algebra.


\begin{prop}
Let $c=(C_{0}, C_{1},\ldots)$ be a 1-truncated conformal algebra. 
Set $C = C_{0} \oplus C_{1} \oplus( \oplus_{i\geq 1} \partial^{i} C_{1}) $ 
then there is a unique structure of a conformal algebra on $C$ 
such that the restrictions of the operations $_{(n)}, n\geq 0$ and $\partial$ to the subspace $C_{\leq 1}$ coincide with the ones given by the 1-truncated conformal algebra structure of $c$.
\end{prop}

\begin{proof}
This is Theorem 9.8 in [GMS1].
\end{proof}


The following proposition explains how to construct the vertex envelope of a conformal algebra, which carries the structure of a vertex algebra. 

\begin{prop} \label{uniofconf}
Let $C$ be an $\mathbb{N}$-graded conformal algebra, 
then the operation 
$[x, y] = \sum_{j\geq 0} \partial^{(j+1)}(x_{(j)} y)$ defines a Lie algebra structure on $C$. 
Moreover there is a unique structure of a vertex algebra on $UC$ 
(the universal enveloping algebra of $C$ with respect to $[,]$) 
such that for any $x\in C$, $z\in TC$ (the tensor algebra of $C$), 
$$ p(xz) = i(x)_{(-1)} p(z) $$ 
where $p: TC\to UC$ is the canonical projection 
and $i :C\to UC$ is the restriction of $p$ to $C\subset TC$.
\end{prop}

\begin{proof}
This is Theorem 8.3 in [GMS1].
\end{proof}

Now we explain how to construct the enveloping algebra $U\mathcal{A}$ of a vertex algebroid $\mathcal{A}$. 
First by Proposition \ref{broidtotrun}, regard $\mathcal{A}$ as a 1-truncated vertex algebra. 
Forget the $_{(-1)}$ operation for a second, it is also a 1-truncated conformal algebra. 
Extend it to an $\mathbb{N}$-graded conformal algebra denoted by $C\mathcal{A}$ 
and construct its vertex envelope $UC\mathcal{A}$. 
But we are not done yet. 
The vertex algebra $UC\mathcal{A}$ is quite big. 
Now the $_{(-1)}$ operation comes into play and it would help us to reduce the size of the vertex algebra enormously. 
This is the next proposition. 


\begin{prop}
Let $\mathcal{A} = (A, T, \Omega, \ldots)$ be a vertex algebroid. 
Denote by $R$ the left $\partial$-ideal of $UC\mathcal{A}$ generated by the following elements
$$ (1_{A} - \mathbf{1} ) u , \quad u\in UC\mathcal{A};$$
$$ a\cdot x - ax, \quad a\in A, \quad x\in A\oplus \Omega; \qquad  a\cdot \tau - a\tau+\gamma(a, \tau), \quad \tau\in T.$$
Then $R$ is a vertex ideal in $UC\mathcal{A}$.  
The vertex algebra $U\mathcal{A}:=UC\mathcal{A}/R$ is called the enveloping algebra of $\mathcal A$ and it satisfies that $U\mathcal{A}_{0} = A$, $U\mathcal{A}_{1} = T\oplus \Omega$. 
\end{prop}

\begin{proof}
See 9.11, 9.12, 9.13 and 9.14 in [GMS1].
\end{proof}

\begin{remark}
$UC\mathcal{A}$ is an associative algebra and a vertex algebra. 
But the algebra multiplication is not induced by the $_{(-1)}$ operation of the vertex algebra structure. Since the left $\partial$-ideal $R$ with respect to the algebra structure is also a vertex ideal, 
the vertex algebra structure of $U\mathcal{A}= UC\mathcal{A}/ R$ 
is inherited from that of $UC\mathcal{A}$, 
but it usually does not inherit the algebra structure of $UC\mathcal{A}$. 
\end{remark}

From now on let us focus on the vertex algebroid 
$\mathcal{A}_{\mf{g}, k} = ( A= \mathcal{R}(G), T= A\otimes \mf{g}, \Omega = \text{Hom}_{A} (T, A), \partial, \gamma, k(,), c) $ 
associated to a simple complex Lie group $G$ 
with Lie algebra $\mf{g}$ and $k\in \mathbb{C}$ (see Proposition \ref{11}). 
The existence of an  embedding of vertex algebras 
$V_{0, k}\hookrightarrow U\mathcal{A}_{\mf{g}, k}$ 
is obvious from the definition of $\mathcal{A}_{\mf{g}, k}$. 
Moreover it is proved in [AG] and as Theorem 2.5 (B. Feigin - E. Frenkel, D.Gaitsgory) in [GMS2] 
that there exists an embedding of vertex algebras 
$V_{0, \bar{k}=-2h^{\vee}-k} \hookrightarrow U\mathcal{A}_{\mf{g}, k}$ 
and the image commutes with $V_{0, k}$ in the appropriate sense. 
In particular 
$U\mathcal{A}_{\mf{g}, k}$ is a $\hat{\mf{g}}\oplus \hat{\mf{g}}$-module of dual levels. 
We prove that for generic $k$, $U\mathcal{A}_{\mf{g}, k}$ 
is isomorphic to the vertex operator algebra $\V$ constructed in Section 2. 

\begin{theorem}\label{main2}
Let $G$ be a simple complex Lie group with Lie algebra $\mf{g}$ and $k \notin \mathbb{Q}$. 
The enveloping algebra $U\mathcal{A}_{\mf{g}, k}$ of the vertex algebroid $\mathcal{A}_{\mf{g}, k}$
decomposes into $\bigoplus_{\lambda\in P^{+}} V_{\lambda, k} \otimes V_{\lambda^{*}, \bar{k}} $ 
as a $\hat{\mf{g}}\oplus \hat{\mf{g}}$-representation. 
Moreover $U\mathcal{A}_{\mf{g}, k}\cong \V$ as vertex operator algebras.
\end{theorem}

The theorem will be proved towards the end of Section 3 after we study the structure of $U\mathcal{A}_{\mf{g}, k}$ carefully by a sequence of lemmas and propositions. First we introduce some notations.

Let $\{\tau_{i}\}$ be an orthonormal basis of $\mf{g}$ 
with respect to  the normalized invariant symmetric bilinear form 
($\tau_{i}$ is the same as $u_{i}$ in Section 2, 
but we use $\tau_{i}$ here to agree with the notations in [GMS1,2]).
Let  $C_{ijk}$ be the structure constants determined by 
$[\tau_{i}, \tau_{j}] = C_{ijk}\tau_{k}$, 
then $C_{ijk} = ([\tau_{i}, \tau_{j}], \tau_{k}) = (\tau_{i}, [\tau_{j}, \tau_{k}]) = (\tau_{j}, [\tau_{k}, \tau_{i}]) = C_{jki} =C_{kij} = -C_{jik}$ 
and $$C_{ipq}C_{jqp} = 2h^{\vee} \delta_{ij}.$$
There are two embeddings of the Lie algebras 
$i_{L}: \mf{g} \hookrightarrow T,  i_{R}: \mf{g} \hookrightarrow T$ 
having left and right invariant vector fields as their images and 
$[i_{L}(x), i_{R}(y)] = 0$ for any $x,y \in \mf{g}$. 
Identify $\mf{g}$ with its image under $i_{L}$ and denote $i_{R}(x)$ by $x^{R}$.  
There exists $a^{ij}\in A$ such that 
$$\tau_{i}^{R} = a^{ij} \tau_{j}$$ 
where $a^{ij} (1) = -\delta_{ij}$, 
i.e. $(a^{ij})$ is the invertible transformation matrix between the two $A$-basis of $T$. 
The following identities are derived in [GMS2]:
$$ \tau_{i} (a^{jk}) = - C_{ipk} a^{jp}, \qquad       \tau_{i}^{R} (a^{jk}) = C_{ijq} a^{qk}.$$

\begin{lemma}
One has 
$a^{ij}a^{ik} = \delta_{jk}$, $a^{ik}a^{jk} = \delta_{ij}$ and $\tau_{i} = a^{ji} \tau_{j}^{R}$.
\end{lemma}

\begin{proof}
For any $s$, we have 
$\tau_{s}^{R} (a^{ij} a^{ik}) = \tau_{s}^{R} (a^{ij}) a^{ik} + a^{ij} \tau_{s}^{R} (a^{ik}) = C_{sip} a^{pj} a^{ik} + C_{siq} a^{ij} a^{qk} = 0 $, 
hence $a^{ij}a^{ik} = \text{ constant } = a^{ij}(1) a^{ik}(1) =( -\delta_{ij}) (-\delta_{ik}) = \delta_{jk}$. Similarly $\tau_{s}(a^{ik}a^{jk}) = -C_{spk}a^{ip}a^{jk} - C_{sqk}a^{ik}a^{jq} =0 $, 
hence $a^{ik} a^{jk} = (-\delta_{ik})(-\delta_{jk})= \delta_{ij}$. 
The third identity follows. 
\end{proof}

Let $\{\omega_{i}\}$ be (left invariant) 1-forms dual to the vector fields $\{\tau_{i}\}$, i.e. $\langle \tau_{i}, \omega_{j} \rangle = \delta_{ij}$. 
$\{\omega_{i}\}$ form an $A$-basis of $\Omega$ and one has according to [GMS2]:
$$\tau_{i} (\omega_{j}) = C_{ijs} \omega_{s}, \qquad   \tau_{i}^{R} (\omega_{j}) =0. $$

$\{\tau_{i}(z)\}$ generate one copy of the $\hat{\mf{g}}$-action 
on $U\mathcal{A}_{\mf{g}, k}$ (of level $k$) 
since ${\tau_{i}}_{(0)}\tau_{j} = C_{ijs}\tau_{s}$ and  ${\tau_{i}}_{(1)}\tau_{j} = k\delta_{ij}$. 
Let 
$$j_{R}( {\tau_{i}}) = \tau_{i}^{R} + \bar{k} a^{ip} \omega_{p}, $$ 
it is shown in [GMS2] that 
$\{ j_{R}(\tau_{i})(z) \}$ generate another copy of the $\hat{\mf{g}}$-action 
on $U\mathcal{A}_{\mf{g}, k}$ of the dual level $\bar{k} = -2 h^{\vee}-k$. 
We use the 'bar' notation to denote it, 
so $\bar{\tau_{i}}(-1) \mathbf{1} = j_{R}(\tau_{i}) = a^{ij}_{(-1)} \tau_{j}(-1)\mathbf{1} + \bar{k} a^{ip}\omega_{p} = \tau_{j}(-1) a^{ij} + \bar{k} a^{ip}\omega_{p}$. 
The second equality is because $\gamma| _{A \times \mf{g}} =0$ and the last equality is because 
$[\tau_{j}(m), a^{ij} _{(n)}]=0$ for any $m,n\in \mathbb{Z}$ since $\tau_{j} (a^{ij}) = -C_{jpj}a^{ip} =0$. 

Fix $\mathcal{A} = \mathcal{A}_{\mf{g}, k} = (A, T, \Omega, \ldots)$.

\begin{lemma} \label{ca}
Let $\{ \widetilde{\omega_{i}} = a^{ij} \omega_{j} \}$ 
be the right invariant 1-forms dual to $\{ \tau_{i}^{R} \}$. 
For any $a, b\in A$, $\omega, \omega' \in \Omega$, $v\in U\mathcal{A}$, we have

\begin{enumerate}
\item 
$ \bar{\tau_{i}}(-1)\mathbf{1} = \tau_{j}(-1) a^{ij} + \bar{k} \widetilde{\omega_{i}}$
\item
$ \tau_{j}(-1)\mathbf{1} = \bar{\tau_{i}}(-1) a^{ij} + k\omega_{j} $
\item 
$[a_{(m)}, b_{(n)}] = [a_{(m)}, \omega_{(n)}] = [\omega_{(m)}, {\omega'}_{(n)}] =0$
\item
$[\tau_{i}(m), a_{(n)}] = (\tau_{i} a)_{(m+n)}$ and 
$ [\bar{\tau_{i}}(m), a_{(n)}] = (\tau_{i}^{R} a)_{(m+n)}$
\item
$[\tau_{i}(m),{\omega_{j}}_{(n)}] = C_{ijs} {\omega_{s}}_{(m+n)} + m\delta_{m+n, 0} \delta_{ij}$
\item
$[\bar{\tau_{i}} (m), \widetilde{\omega_{j}}_{(n)}] = C_{ijs} \widetilde{\omega_{s}}_{(m+n)} + m\delta_{m+n, 0}\delta_{ij}$
\item
$[\bar{\tau_{i}}(m), {\omega_{j}}_{(n)}] = m a^{ij}_{(m+n-1)}$  and 
$[\tau_{i}(m), \widetilde{\omega_{j}}_{(n)}] = ma^{ji}_{(m+n-1)}$
\item 
$\bar{\tau_{i}} (m) v  = \sum_{l\geq 0} \{ \tau_{j}(-1-l) a^{ij}_{(m+l)} v + a^{ij}_{(m-1-l)} \tau_{j} (l) v \}
+\bar{k} \widetilde{\omega_{i}}_{(m)}v $
\item
$\tau_{j}(m)v = \sum_{l\geq 0} \{\bar{\tau_{i}} (-1-l) a^{ij}_{(m+l)} v + a^{ij}_{(m-1-l)} \bar{\tau_{i}}(l) v \}+ k {\omega_{j}}_{(m)} v$
 \end{enumerate}
 \end{lemma}

\begin{proof}
Everything follows from either the commutator or iterate formula once the second identity is verified which is straightforward computation. 
Since 
$\partial a^{ij}=\tau_{k}(a^{ij})\omega_{k}=-C_{kpj}a^{ip}\omega_{k}$,  
by the iterate formula we have 
$\bar{\tau_{i}}(-1) a^{ij} 
= \tau_{k}(-1) a^{ik}_{(-1)} a^{ij} + a^{ik}_{(-2)} (\tau_{k} a^{ij}) +\bar{k} a^{ip}\omega_{p} a^{ij} 
= \tau_{k}(-1)\delta_{kj} + (-C_{lpk}a^{ip}\omega_{l}) \\ (-C_{kqj}a^{iq}) + \bar{k}\omega_{p}\delta_{pj} 
= \tau_{j}(-1)\mathbf{1} + C_{kqj}C_{lpk}\delta_{pq}\omega_{l}+\bar{k}\omega_{j} 
= \tau_{j}(-1)\mathbf{1} +2h^{\vee}\delta_{jl}\omega_{l}+\bar{k}\omega_{j} 
= \tau_{j}(-1)\mathbf{1}-k\omega_{j}$, 
hence $\tau_{j}(-1)\mathbf{1}=\bar{\tau_{i}}(-1)a^{ij}+k\omega_{j}$. 
\end{proof}

Let $B'$ be the subalgebra generated by $A$ and $\partial^{(i)}\Omega, i\geq 0$ in $UC\mathcal{A}$. 
In fact $B' = \text{Sym}_{\mathbb{C}}(A\oplus (\oplus_{i\geq 0}\partial^{(i)} \Omega))$. 
Let $B$ be the image of $B'$ in $U\mathcal{A}$, 
as we'll see later that 
it plays an important role in our proof of Theorem \ref{main2}.

As an abelian vertex subalgebra of $U\mathcal{A}$, 
$B$ is the jet algebra associated to the abelian vertex algebroid 
$(A, \Omega, \partial)$ (cf. [GMS1] 9.3). 
$B$ is a free commutative unital $A$-algebra 
with generators  $\{ \partial^{(j)}\omega_{i}, j\geq 0, i\}$ 
(or $\{ \partial^{(j)}\widetilde{\omega_{i}}, j\geq 0, i\}$) and 
is $\mathbb{N}$-graded with 
$B_{0}=A$, 
$B_{1}=\Omega=\oplus_{i} A\omega_{i}$( or $\oplus_{i} A\widetilde{\omega_{i}}$), 
$B_{2}=(\oplus_{i} A\partial \omega_{i})\oplus  (\oplus_{i,j} A\omega_{i}\omega_{j})$ 
(or $(\oplus_{i} A \partial \widetilde{\omega_{i}})\oplus (\oplus_{i,j} A\widetilde{\omega_{i}}\widetilde{\omega_{j}})$) 
and etc. 

\begin{lemma}\label{AB}
$B$ is the vertex subalgebra of $U\mathcal{A}$ generated by $U\mathcal{A}_{0}=A$.
\end{lemma}

\begin{proof}
Since 
$a^{ij}_{(-2)}a^{il}
= (-C_{kpj}a^{ip}\omega_k) a^{il}
= C_{kjl}\omega_{k}$, 
we deduce that 
$\omega_{r}=\frac{1}{2h^{\vee}} C_{rlj} a^{ij}_{(-2)} a^{il}$.
The lemma follows. 
\end{proof}

Let $U(\hat{\mf{g}}, k) = U(\hat{\mf{g}})/(\underline c - k )$. 

\begin{prop} \label{PBW}
$B$ is closed under the actions of $\hat{\mf{g}}_{\geq 0}$ and $\bar{\hat{\mf{g}}}_{\geq 0}$. 
As a $U(\hat{\mf{g}}, k)$-module, 
$U\mathcal{A} = U(\hat{\mf g}, k) \otimes_{U(\hat{\mf g}_{\geq 0}, k)} B$; 
as a $U(\bar{\hat{\mf{g}}}, \bar{k})$-module, 
$U\mathcal{A} = U(\bar{\hat{\mf g}}, \bar k) \otimes_{U(\bar{\hat{\mf g}}_{\geq 0}, \bar k)} B$. 
In particular as a vector space 
$U\mathcal{A} = U(\hat{\mf{g}}_{-}) \otimes B = U(\bar{\hat{\mf{g}}}_{-}) \otimes B$. 
Furthermore $B^{\hat{\mf{g}}_{>0}} = B^{ \bar{\hat{\mf{g}}}_{>0}} =A$, 
$U\mathcal{A}^{\hat{\mf{g}}_{>0}} = U( \bar{\hat{\mf{g}}}_{-}) \otimes A$ and 
$U\mathcal{A}^{ \bar{\hat{\mf{g}}}_{>0}} = U(\hat{\mf{g}}_{-}) \otimes A$.
\end{prop}

\begin{proof}
It follows from the commutation relations in Lemma \ref{ca} that 
$B$ is closed under the $\hat{\mf g}_{\geq 0}$- and $\bar { \hat{\mf g} }_{\geq 0}$-actions. 
Write the vectors in $B$ as polynomials in $\{ \partial^{(j)} \omega_i\}$ with coefficients in $A$, 
we say that the highest degree of 
$b = a\partial^{(p_{1})} \omega_{i_{1}} \cdots \partial^{(p_{l})}\omega_{i_{l}}$ 
is $p_{1}+1$ if $p_{1}\geq p_{2}\geq \cdots \geq p_{l}$. 
Then $\tau_{i}(m), m>0$ acts as $m$ times differentiation with respect to the term 
$\partial^{(m-1)} \omega_{i}$ on any vector of highest degree no more than $m$. 
In particular $\tau_{i}(m)$ kills all vectors of highest degree strictly less than $m$. 
Based on this observation it is not difficult to see that $B^{\hat{\mf{g}}_{>0}} = A$. 
For example suppose that 
$b_{2}=\sum a_{i}\partial \omega_{i} +\sum_{i\leq j} a_{ij} \omega_{i}\omega_{j}\in B_{2}$ 
is killed by $\hat{\mf{g}}_{>0}$ where $a_{i}, a_{ij}\in A$,  
then $a_{i}= \frac{1}{2}\tau_{i}(2) b_{2}=0$, 
therefore $b_2 = \sum_{i\leq j} a_{ij} \omega_{i}\omega_{j}$, 
but again $a_{ii}=\frac{1}{2}\tau_{i}(1)\tau_{i}(1)b_{2} =0$ and 
$a_{ij} = \tau_{i}(1)\tau_{j}(1) b_{2} =0$ for $i<j$, 
hence $b_{2}=0$.  
Similarly we can show that $B^{\bar {\hat{ \mf g}}_{>0}} = A$, 
but when dealing with the $\bar{\hat{\mf{g}}}_{\geq0}$-action on $B$, 
it's more natural to use another set of $A$-generators $\{\partial^{(j)}\widetilde{\omega_{i}}\}$ of $B$ 
since it will prevent the appearance of structure functions $a^{ij}$ 
(see (6) (7) of Lemma \ref{ca}). 
It follows from Theorem 9.16 in [GMS1] that as a vector space 
$U\mathcal{A} \cong U(\hat{\mf{g}}_{<0}) \otimes B \cong U(\bar{\hat{\mf{g}}}_{<0}) \otimes B$. 
Roughly speaking $\{\tau_{i}\}$ (or $\{ \tau_i^R \}$) form an $A$-basis of $T$ 
and the function part can always be moved to the right and absorbed into $B$.  
This implies that as both a $\hat{\mf{g}}$- and $\bar{\hat{\mf{g}}}$-module, 
$U\mathcal{A}$ is induced from $B$. 
Since $U\mathcal{A} = U(\bar{\hat{\mf{g}}}_{-}) \otimes B$ (resp. $ U(\hat{\mf g}_-) \otimes B$) 
and the two copies of the $\hat{\mf{g}}$-actions commute with each other, 
it follows from $B^{\hat{\mf{g}}_{>0}} = A$ (resp.  $B^{\bar {\hat{\mf{g}}}_{>0}} = A$) 
that $U\mathcal{A}^{\hat{\mf{g}}_{>0}} = U( \bar{\hat{\mf{g}}}_{-}) \otimes A$
(resp. $U\mathcal{A}^{ \bar{\hat{\mf{g}}}_{>0}} = U(\hat{\mf{g}}_{-}) \otimes A$).  
\end{proof}

\begin{corollary}
Let $W$ be the vertex algebra generated by the quantum fields 
$a(z) = \sum a_{(n)} z^{-n-1}, a\in A$ of conformal weight $0$ 
and $\tau_{i}(z) = \sum {\tau_i}_{(n)} z^{-n-1}$, 
$\omega_{i}(z) = \sum {\omega_i}_{(n)} z^{-n-1}$ of conformal weight $1$
with the OPE: 
$$\tau_{i}(z_{1}) a(z_{2}) \sim \frac{\tau_{i}a(z_{2})}{z_{1}-z_{2}}$$
$$\tau_{i}(z_{1})\tau_{j}(z_{2})\sim \frac{C_{ijs}\tau_{s}(z_{2})}{z_{1}-z_{2}} + \frac{k\delta_{ij}}{(z_{1}-z_{2})^{2}}$$
$$\tau_{i}(z_{1})\omega_{j}(z_{2})\sim \frac{C_{ijs}\omega_{s}(z_{2})}{z_{1}-z_{2}} + \frac{\delta_{ij}}{(z_{1}-z_{2})^{2}}$$ 
$$a(z_{1})b(z_{2})\sim 0 \qquad a(z_{1})\omega_{j}(z_{2})\sim 0 \qquad \omega_{i}(z_{1})\omega_{j}(z_{2})\sim 0$$
and let $I$ be the vertex $\partial$-ideal generated by
$a_{(-1)}a'_{(-1)} \mathbf 1 - (a a')_{(-1)} \mathbf 1$ and 
$a_{(-2)} \mathbf 1 - (\tau_i a)_{(-1)} {\omega_i}_{(-1)} \mathbf 1$, $a, a' \in A$, 
then $W/I \cong U\mathcal{A}$. 
\end{corollary}

\begin{proof}
It follows from Lemma \ref{ca} and Proposition \ref{PBW}.
\end{proof}

\begin{remark}
We can certainly formulate other descriptions. 
For example if we replace $\omega_{j}(z)$ by $\widetilde{\omega_{j}}(z)$, 
then all the OPE remain the same except the one involving 
$\tau_{i}(z_{1})$ and  $\omega_{j}(z_{2})$ 
which ought to be replaced by 
$\tau_{i}(z_{1})\widetilde{\omega_{j}}(z_{2}) \sim \frac{a^{ji}(z_{2})}{(z_{1}-z_{2})^{2}}$, 
and part of the generators of the ideal should be replaced by 
$a_{(-2)} \mathbf 1 - (\tau_j^R a)_{(-1)} { \widetilde{\omega_j} }_{(-1)} \mathbf 1$. 
Two other options are to use $\bar{\tau_{i}}(z)$ and $\omega_{i}(z)$ 
or $\bar{\tau_{i}}(z)$ and $\widetilde {\omega_i}(z)$. 
\end{remark}


Set 
$$\varpi = \tau_j (-1)\omega_j + \frac{\bar k }{2}  \omega_s^2 \in U\mathcal A_2, $$
we will show that $\varpi$ is the \emph{Virasoro element} in $U\mathcal{A}$ 
(for generic $k$, $\varpi$ is the same as $\omega\in \V_{2}$ in Section 2). 

\begin{lemma} 
One has 
$a^{ij}_{(-3)} a^{ij} = h^{\vee} \omega_{s}^{2}$ and $\omega_{s}^{2} = \widetilde{\omega_{t}}^{2}$.
\end{lemma}

\begin{proof}
Since $a^{ij}_{(-2)}a^{ij}=0$, 
we have $0=\partial(a^{ij}_{(-2)}a^{ij})= 2 a^{ij}_{(-3)}a^{ij} + (\partial a^{ij}) (\partial a^{ij})$. 
Since $ (\partial a^{ij}) (\partial a^{ij}) = (\tau_{k} a^{ij}) (\tau_{l} a^{ij}) \omega_{k}\omega_{l} = C_{kpj} a^{ip} C_{lqj} a^{iq} \omega_{k}\omega_{l} = -2 h^{\vee} \omega_{k}^{2}$, 
we deduce that $a^{ij}_{(-3)} a^{ij} = h^{\vee} \omega_{k}^{2}$. 
The other identity is easily checked. 
\end{proof}

\begin{prop} 
One has 
$$\varpi = \bar{\tau_{i}}(-1) \widetilde{\omega_{i}} + \frac{k}{2} \omega_s^2. $$
If $\varkappa = k+ h^\vee \neq 0$, then
$$\varpi= \frac{1}{2\varkappa} \sum_i \tau_{i}(-1)\tau_{i}(-1)\mathbf{1} - \frac{1}{2\varkappa} \sum_j \bar{\tau_{j}}(-1) \bar{\tau_{j}}(-1)\mathbf{1}. $$
\end{prop}

\begin{proof}
By (9) of Lemma \ref{ca}, we have 
$\tau_{j}(-1) \omega_{j} = \bar{\tau_{i}}(-1) a^{ij} \omega_{j} + a^{ij}_{(-2)} \bar{\tau_{i}}(0) \omega_{j} + a^{ij}_{(-3)} \bar{\tau_{i}}(1) \omega_{j} + k\omega_{j}^{2}= \bar{\tau_{i}}(-1) \widetilde{\omega_{i}} 
+ 0 + h^{\vee} w_{s}^{2} + k\omega_{j}^{2} =  \bar{\tau_{i}}(-1) \widetilde{\omega_{i}} + \varkappa \omega_{s}^{2}$, 
hence $\varpi =  \bar{\tau_{i}}(-1) \widetilde{\omega_{i}} + \varkappa \omega_{s}^{2} - \frac{\varkappa + h^{\vee}}{2} \omega_{s}^{2} =  \bar{\tau_{i}}(-1) \widetilde{\omega_{i}} + \frac{k}{2} \omega_{s}^{2}$.
Since 
$\sum \tau_{j}(-1)\tau_{j}(-1)\mathbf{1} - \sum \bar{\tau_{i}}(-1) \bar{\tau_{i}}(-1)\mathbf{1} = \sum_{j} \tau_{j}(-1)(\tau_{j}(-1)\mathbf{1} -\sum_{i} \bar{\tau_{i}}(-1) a^{ij}) 
- \sum_{i} \bar{\tau_{i}}(-1)(\bar{\tau_{i}}(-1)\mathbf{1} - \sum_{j} \tau_{j}(-1) a^{ij}) =  k \sum \tau_{j}(-1) \omega_{j} - \bar{k} \sum \bar{\tau_{i}}(-1)\widetilde{\omega_{i}} = k (\varpi - \frac{\bar{k}}{2} \omega_{s}^{2}) - \bar{k} (\varpi- \frac{k}{2}\omega_{s}^{2}) = 2 \varkappa \varpi$, 
the other identity follows. 
\end{proof}

\begin{lemma}\label{lem1}
One has the following: 
\begin{enumerate}
\item 
$a_{(0)}\varpi = -\partial a$ 
and  $a_{(n)}\varpi = 0$ for any $a\in A$, $n\geq 1$; 
\item 
$\tau_{i}(0) \varpi =0$, 
$\tau_{i}(1)\varpi = \tau_{i}(-1)\mathbf{1}$ and 
$\tau_{i}(n) \varpi =0$ for all $n\geq 2$; 
\item 
$\bar{\tau_{i}}(0) \varpi =0$, 
$\bar{\tau_{i}}(1) \varpi = \bar{\tau_{i}}(-1)\mathbf{1}$ and 
$\bar{\tau_{i}}(n)\varpi =0$ for all $n\geq 2$; 
\item 
${\omega_{i}}_{(0)} \varpi =0$, 
${\omega_{i}}_{(1)}\varpi = \omega_{i}$ and 
${\omega_{i}}_{(n)}\varpi =0$ for all $n\geq 2$. 
\end{enumerate}
\end{lemma}

\begin{proof}
Straightforward computation. It follows from Lemma \ref{ca} and the definition of $\varpi$.
\end{proof}

\begin{prop}\label{Vira1}
Let $Y(\varpi, z) =\sum_{n\in \mathbb{Z}} \mathcal{L}_{n} z^{-n-2}$, 
then for any $a\in A, x, y \in \mf g, \omega \in \Omega$ and $m, n\in \mathbb Z$, one has the following:
\begin{enumerate}
\item 
$[\mathcal{L}_{m}, a_{(n)}] = - (m+n+1) a_{(m+n)}$; 
\item 
$[\mathcal{L}_{m},  x(n)] = -n x(m+n)$; 
\item 
$[\mathcal{L}_{m}, \bar y (n)] = -n \bar y(m+n)$; 
\item 
$[\mathcal{L}_{m}, {\omega}_{(n)}] = -n {\omega}_{(m+n)}$. 
\end{enumerate}
\end{prop}

\begin{proof}
It follows from Lemma \ref{lem1}  and Borcherds' commutator formula.
\end{proof}

\begin{lemma}\label{Vira2}
One has 
$\varpi_{(0)}\varpi = \partial \varpi$, 
$\varpi_{(1)} \varpi = 2 \varpi$, 
$ \varpi_{(2)}\varpi =0$ and 
$\varpi_{(3)} \varpi = \text{dim} \mf{g} $.
\end{lemma}

\begin{proof}
Straightforward computation using Lemma \ref{lem1} or Proposition \ref{Vira1}.
\end{proof}

\begin{prop}
One has 
$[\mathcal{L}_{m}, \mathcal{L}_{n}] = (m-n) \mathcal{L}_{m+n} + \frac{m^{3}-m}{12} \delta_{m+n, 0} 2\,\text{dim}\mf{g}$, 
$\mathcal{L}_{-1} = \partial$ and that 
$\mathcal{L}_{0}$ acts as the gradation operator.  
Hence $\varpi$ is a Virasoro element and 
$U\mathcal{A}$ is a vertex operator algebra of rank $2\, \text{dim} \mf{g}$. 
\end{prop}

\begin{proof}
The Virasoro commutation relations follow from Lemma \ref{Vira2}. 
The fact that $\mathcal{L}_{-1} = \partial$ can be first checked on the generators 
$a\in A$, $\tau_{i}$ and $\omega_{i}$, then on the whole vertex algebra $U\mathcal A$. 
The commutation relations for $\mathcal{L}_{0}$ in Proposition \ref{Vira1}  
imply that $\mathcal{L}_{0}$ acts on $U\mathcal{A}$ as the gradation operator. 
\end{proof}


\begin{define} \label{bfonA}
Define a symmetric bilinear form $(, )$ on $A$ as follows: for $u\otimes u^{*} \in V_{\lambda}\otimes V_{\lambda}^{*}$ and $v^{*}\otimes v\in V_{\mu}^{*}\otimes V_{\mu}$, set $( u\otimes u^{*}, v^{*}\otimes v ) = \frac{1}{\textrm{dim $V_{\lambda}$}}( u, v^{*}) ( u^{*}, v)$ if $\lambda =\mu$; $0$ otherwise, then extend it bilinearly to $A\times A \to \mathbb{C}$. In particular $( 1,1) = 1$.
\end{define}

\begin{lemma} \label{bfA}
The bilinear form $(, )$ defined above satisfies the following: 
\begin{enumerate}
\item $(, )$ is non-degenerate; 
\item $( a, b)  = ( 1, ab )$, in particular $( ab, c ) = ( b, ac ) $ for any $a, b, c\in A$; 
\item $( x\cdot a, b) = -( a, x \cdot b )$ and $( y^R \cdot a, b) = -( a, y^R \cdot b)$
for any $x, y\in \mf{g}$. 
\end{enumerate}
\end{lemma}

\begin{proof}
Non-degeneracy is obvious. 
It suffices to prove (2) for $a = u\otimes u^{*} \in V_{\lambda}\otimes V_{\lambda}^{*}$ 
and $b = v\otimes v^{*}\in V_{\mu}\otimes V_{\mu}^{*}$. 
Note that 
$V_{\lambda \mu}^{0} \cong (V_{\lambda}^{*} \otimes V_{\mu}^{*} )^{\mf{g}} 
\cong \text{Hom}_{\mf{g}}(V_{\lambda}, V_{\mu}^{*}) 
= \mathbb{C} \text{  if $\mu = \lambda^{*} (= \omega_0 \lambda) $}; 
0 \text{  otherwise}$. 
Hence if $\mu\neq \lambda^{*}$, 
then $ab$ has no $0$-component, in which case $(1, ab)=0= (a, b)$. 
If  $\mu = \lambda^{*}$, then $V_{\lambda \lambda^{*}}^{0} = \mathbb{C}$ 
and it is generated by the evaluation map $f: w \otimes w^{*} \mapsto (w, w^{*})$. 
It's not difficult to check that the evaluation map 
$g\in V_{\lambda^{*} \lambda}^{0}$ has dual 
$g^{*} : \mathbb{C} \to V_{\lambda}\otimes V_{\lambda}^{*}; 1 \mapsto \sum_{i} w_{i} \otimes w_{i}^{*}$, 
where $\{w_i \}$ and $\{ w_i^*\}$ are dual basis of $V_{\lambda}$ and $V_{\lambda}^{*}$. 
Hence $f g^{*} = \sum_{i} (w_{i}, w_{i}^{*} ) = \text{dim} V_{\lambda}$, 
which means that $(f, g) = \text{dim} V_{\lambda}$. 
By Proposition \ref{peterweyl} the $0$-component of $ab$ is 
$\frac{1}{\text{dim $V_{\lambda}$}} f(u \otimes v) \otimes g( u^{*} \otimes v^{*}) = \frac{1}{\text{dim $V_{\lambda}$}} (u, v) (u^{*}, v^{*}) $, 
therefore $(1, ab)= (1, (ab)_{0}) =(a, b)$. (3) is obvious.
\end{proof}

\begin{remark}
If we regard $A$ as the space of representative functions on the maximal compact Lie subgroup of $G$,
then $(a, b)$ is nothing else but the integral of $ab$ with respect to the Haar measure. 
\end{remark}

Let $\sigma$ be the anti-involution of $\hat{\mf{g}}$ defined by 
$x(m) \mapsto -x(-m), \underline c \mapsto \underline c$. 
Recall $U(\hat{\mf{g}}, k) = U(\hat{\mf{g}}) /(\underline c = k)$, 
then $\sigma$ can be extended to an antiautomorphism of $U(\hat{\mf{g}}, k)$ 
which we still denote by $\sigma$.  
We also have the PBW decomposition 
$U(\hat{\mf{g}}, k) \cong U(\hat{\mf{g}}_{-}) \otimes U(\mf{g}) \otimes U(\hat{\mf{g}}_{+})$. 

Let's recall the contragredient module $V^{c}$ of a vertex operator algebra $V$ from [FHL]. 
As a vector space, $V^{c} = \oplus_n V_{n}^{*}$ is the restricted dual of $V$. 
The $V$-module structure is given by 
$(Y(v, z) u^{*}, u) = (u^{*}, Y(e^{z\mathcal{L}_{1}} (-z^{-2})^{\mathcal{L}_{0}} v, z^{-1}) u)$ 
for any $u^{*}\in V^{c}, u, v\in V$, which is equivalent to 
$( v_{( 2n-m-2)} u^{*}, u ) = (-1)^{n}$
$( u^{*}, \sum_{k\geq 0} \frac{1}{k!} (\mathcal{L}_{1}^{k} v)_{(m-k)} u) $ 
for homogeneous vector $v$ of degree $n$ and any $m\in \mathbb{Z}$. 
In particular for $a\in V_{0}$, 
one has $(a_{(m)} u^{*}, u) = (u^{*}, a_{(-m-2)} u )$ and for $v \in V_{1}$, 
one has $(v_{(m)}u^{*}, u ) = - (u^{*}, \{v_{(-m)} + (\mathcal{L}_{1}v)_{(-m-1)}\} u)$.

\begin{prop}\label{contragrediant}
There exists a unique non-degenerate symmetric bilinear form 
$(, )$ on $U\mathcal{A}$ such that it satisfies the following: 
\begin{enumerate}
\item 
it is an extension of the form $(, )$ on $A = U\mathcal{A}_{0}$; 
\item 
$( U\mathcal{A}_{m}, U\mathcal{A}_{n}) =0$ except for $m=n \geq 0$; 
\item 
$(Y(v, z) v_1, v_2 ) = ( v_1, Y( e^{ z\mathcal{L}_{1}} (-z^{-2})^{\mathcal{L}_{0}}v, z^{-1}) v_2 )$ 
for any $v, v_1, v_2 \in U\mathcal{A}$, 
in particular 
$a_{(n)}^{*} = a_{(-n-2)}$, 
$\omega_{(n)}^* = - \omega_{(-n)}$, 
$x(n)^* = -x(-n)$, 
$\bar y(n)^* = - \bar y(-n)$ 
for any $a\in A, \omega\in \Omega, x, y\in \mf g$. 
\end{enumerate}
\end{prop}

\begin{proof}
Let's construct the bilinear form $(, )$ explicitly. 
First we extend the form $(, )$ in Definition \ref{bfonA} to $B$ as follows: 
$(B_{m}, B_{n}) = 0$ except when $m=n=0$, 
in which case it coincides with the form $(, )$ on $A$. 
By Proposition \ref{PBW}, $U\mathcal{A} = U(\hat{\mf{g}}_{-}) \otimes B$. 
Given any $P, Q\in U(\hat{\mf{g}}_{-})$ and $b, \tilde{b}\in B$, 
write $\sigma(P)Q = \sum_i R_i^- R_i^0 R_i^+$ for some 
$R_i^- \in U(\hat{\mf{g}}_{-})$, $R_i^0 \in U(\mf g)$ and $R_i^+ \in U(\hat{\mf{g}}_{+})$, 
then define $(Pb, Q\tilde{b}) = \sum_i (\sigma(R_i^-)b, R_i^0 R_i^+ \tilde{b})$. 
Note that the right hand side is well defined since  
both $\sigma(R_i^-) b$ and $R_i^0 R_i^+  \tilde{b}$ are in $B$. 
Since 
$\sigma(Q)P =\sigma(\sigma(P)Q) = \sum_i \sigma(R_i^- R_i^0 R_i^+) 
= \sum_i \sigma(R_i^+) \sigma(R_i^0) \sigma(R_i^-)$, 
by definition 
$(Q\tilde{b}, Pb)= \sum_i ( R_i^+ \tilde{b}, \sigma(R_i^0) \sigma(R_i^-) b) $. 
Then it follows from (3) of Lemma \ref{bfA} that 
$( Pb, Q\tilde{b} ) = ( Q\tilde{b}, Pb )$, 
i.e. the bilinear form $(, )$ we defined is in fact symmetric. 
(2) is obvious. 

From the construction it is easy to check that $x(n)^* = - x(-n)$. 
To show that $a_{(n)}^{*} = a_{(-n-2)}$, i.e.
$(a_{(n)} v_1, v_2) =(v_1, a_{(-n-2)} v_2)$ for any $v_1, v_2\in U\mathcal A$, 
first note that it is true for any $v_1 =b \in B$ and $v_2=\tilde b \in B$
because if $n=-1$, it follows from (2) of Lemma \ref{bfA}; 
otherwise assume $n\geq 0$ without loss of generality, 
then $a_{(n)} b=0$ and $(b, a_{(-n-2)} \tilde{b})=0$ 
since $a_{(-n-2)} \tilde{b} \in \oplus_{n\geq 1}B_{n}$. 
Next we prove that it is true for any $v_1 = b\in B$ and $v_2\in U\mathcal A$. 
Suppose it holds for any $b\in B$ and $v_2$ of length $\leq l$, 
then for any $x\in \mf g, m>0$, we have 
$(a_{(n)}b, x(-m)w)= -(x(m)a_{(n)}b,w) = -((x\cdot a)_{(m+n)}b, w) - (a_{(n)}x(m)b, w)$. 
Since $x(m)b\in B$, by induction it equals 
$- (b, (x\cdot a)_{(-2-m-n)} w)- (x(m)b, a_{(-n-2)} w)
= -(b, (x\cdot a)_{(-2-m-n)} w) +(b, x(-m)a_{(-n-2)} w) 
= (b, a_{(-n-2)} x(-m)w) $, 
which means that it holds for any $b\in B$ and $v_2$ of length $\leq l+1$. 
Finally by similar induction on the length of $v_1$, 
we can show that it is true for any $v_1, v_2\in U\mathcal A$. 
Since 
$\mathcal L_1 x(-1)\mathbf 1 =0$, 
we actually proved (3) for $v = a \in A$ and $x(-1) \mathbf 1$. 
Since they generate the whole vertex operator algebra,  
(3) holds for any $v\in U\mathcal A$. 
Set $v = \omega $ or $\bar y(-1)\mathbf 1$, 
we obtain 
$\omega_{(n)}^* = - \omega_{(-n)}$  and 
$\bar y(n)^*=-\bar y(-n)$
because 
$\mathcal L_1 \omega = \mathcal{L}_{1}\bar y(-1)\mathbf{1} =0$.

Observe that $(,)$ is uniquely determined by (1) and (3), 
so it only remains to prove the non-degeneracy. 
By (2) it suffices to prove it on each fixed level. 
On the top level $U\mathcal{A}_{0}=A$, it is true by (1) of Lemma \ref{bfA}. 
Write $U(\hat{\mf{g}}_{-}) = \oplus_{n\geq 0} U(\hat{\mf{g}}_{-}) [ n]$, 
then $U\mathcal{A}_{1} = U(\hat{\mf{g}}_{-})[1] A \oplus B_{1}$. 
It helps to also decompose it into $U(\bar{\hat{\mf{g}}}_{-})[1] A\oplus B_{1}$. 
Since the two affine actions commute, 
the only nontrivial pairings on $U\mathcal{A}_{1}$ 
are between $U(\hat{\mf{g}}_{-})[1] A$ and $B_{1}$, 
and between $B_{1}$ and $U(\bar{\hat{\mf{g}}}_{-})[1] A$.  
Suppose $v=\sum \tau_i(-1) a_i +\sum c_j \widetilde{\omega_j} \in \text{Ker}(,)|_{U\mathcal A_1}$ for some $a_i, c_j\in A$, 
then for any $s$ and $c\in A$, we have 
$0 = (v, c\omega_{s}) = -\sum (a_{i}, \tau_{i}(1) (c\omega_{s})) 
= - \sum(a_{i}, \delta_{is} c)= -(a_{s}, c)$, 
hence $a_{s}=0$ for any $s$. 
Moreover for any $t$ and $a\in A$, we have 
$0 = (v, \bar{\tau_{t}}(-1) a) = -\sum (\bar{\tau_{t}}(1) (c_{j}\widetilde{\omega_{j}}), a)
=- \sum(c_{j}\delta_{jt}, a)=-(c_{t}, a)$, 
hence $c_{t}=0$ for any $t$. 
Therefore $v=0$, i.e. $\text{Ker}(,)|_{U\mathcal{A}_{1}}=0$. 

Generally, choose an ordered basis of $\hat{\mf{g}}_{-}$: 
$\mathcal X =\{\cdots > \tau_{1}(-n) > \tau_{2}(-n) > \cdots > \tau_{\textrm{dim}\mf{g}}(-n) > \tau_{1}(-n+1) > \tau_{2}(-n+1) > \cdots > \tau_{\textrm{dim}\mf{g}}(-n+1) > \cdots > \tau_{1}(-1) > \tau_{2}(-1) >\cdots >\tau_{\text{dim}\mf{g}}(-1)\}$. 
Then $\mathcal Y = \{x_{1}x_{2}\cdots x_{r} : x_{i}\in \mathcal X, x_1 \geq \cdots \geq x_{i}\geq x_{i+1}\geq \cdots \geq x_r \}$ 
form an encyclopedically ordered basis of $U(\hat{\mf{g}}_{-})$. 
To any $y\in \mathcal{Y}$, 
associate a $\vartheta y\in B$ replacing $\tau_{i}(-n)$ by $\partial^{(n-1)} \omega_{i}$, 
for example 
$\vartheta (\tau_{1}(-2)^{2}\tau_{2}(-2)\tau_{1}(-1)\tau_{2}(-1)^{3}) 
=( \partial \omega_{1})^{2} (\partial \omega_{2}) \omega_{1}\omega_{2}^{3}$. 
By the proof of Proposition \ref{PBW}, it is not difficult to see that 
$\sigma(y_{2}) \vartheta(y_{1})=0$ if $y_{1} < y_{2} \in \mathcal Y$ and 
$\sigma(y_{1}) \vartheta(y_{1})$ equals a nonzero constant. 
For a fixed $n >0$, decompose $U\mathcal A_n$ in two different ways: 
$U\mathcal{A}_{n} = \oplus_{0\leq l\leq n} U(\hat{\mf{g}}_{-})[l]\otimes B_{n-l}
= \oplus _{0\leq l\leq n} U(\bar{\hat{\mf{g}}}_{-})[l] \otimes B_{n-l}$. 
The only nontrivial pairings are between $U(\hat{\mf{g}}_{-})[l] \otimes B_{n-l}$ 
and $ U(\bar{\hat{\mf{g}}}_{-})[n-l]\otimes B_l$ for various $l$. 
Suppose $v= \sum_{0\leq l\leq n} \sum_{1\leq i \leq m_l}  y_i^l b_{n-l}^{i} \\ 
\in \text{Ker}(,)|_{U\mathcal A_n}$ 
where $y_1^l< y_2^l< \cdots \in \mathcal Y \cap U(\hat{\mf{g}}_{-})[l]$ and 
$b_l^i \in B_l$, 
then for a fixed $l$ and any $a\in A, y\in U(\bar{\hat{\mf{g}}}_{-})[n-l]$, we have 
$0 = (v, y a\vartheta (y_1^l)) 
=(\sum_i y_i^l b_{n-l}^i, y a\vartheta( y^{l}_{1})) 
= \sum_i (b_{n-l}^{i}, y \sigma(y_i^l) a \vartheta( y^{l}_{1})
= C (b_{n-l}^{1}, y a) = C (\sigma(y) b_{n-l}^{1}, a)$ 
where $C$ is a nonzero constant. 
Since $a$ is arbitrary, we have 
$\sigma(y) b_{n-l}^{1} =0$ (for any $y \in U(\bar {\hat {\mf g}}_-) [n-l]$). 
By the proof of Proposition \ref{PBW}, this implies that 
$b_{n-l}^1 =0$. 
Inductively we can show that $b_{n-l}^i = 0$ for all $i$, 
hence $v = 0$, i.e. $(,)|_{U\mathcal A_n}$ is non-degenerate.
\end{proof}

Set $A_{\lambda}= V_{\lambda}\otimes V_{\lambda}^{*}$. 
If we take $\{\partial^{(j)} \omega_{i}\}$ as the (free) $A$-generators of $B$, 
then each $b\in B$ can be written as polynomials in these generators with coefficients in $A$. 
Let $B_{\lambda}$ be the collection of those with coefficients in $A_{\lambda}$. 
Similarly let $\widetilde{B_{\lambda}}$ be all the polynomials in generators 
$\partial^{(j)}\widetilde{\omega_{i}}$ with coefficients in $A_{\lambda}$, 
then as a vector space 
$B=\oplus_{\lambda} B_{\lambda} = \oplus_{\lambda} \widetilde{B_{\lambda}}$. 
If we decompose $U\mathcal{A} = \oplus_{\lambda} U(\hat{\mf{g}}_{-})\otimes B_{\lambda}$, 
since $B_{\lambda}$ is closed under the $\hat{\mf{g}}_{\geq 0}$-action, 
it is not difficult to see that 
$(U(\hat{\mf{g}}_{-})\otimes B_{\lambda}, U(\hat{\mf{g}}_{-})\otimes B_{\mu}) = 0$ 
except for $\mu =\lambda^{*}$, 
in addition 
$(U(\hat{\mf{g}}_{-})\otimes B_{\lambda})_{n}$ now becomes finite dimensional. 
The same is true for 
$(U(\hat{\mf{g}}_{-}) \otimes \widetilde{B_{\lambda}}, U(\bar{\hat{\mf{g}}}_{-})\otimes B_{\mu})$ 
and 
$(U(\bar{\hat{\mf{g}}}_{-})\otimes \widetilde{B_{\lambda}}, U(\bar{\hat{\mf{g}}}_{-})\otimes \widetilde{B_{\mu}})$. 

We can make the following identifications: 
 $ U(\hat{\mf{g}}, k) \otimes_{U(\hat{\mf{g}}_{\geq 0}, k)} A_{\lambda} 
 = U(\hat{\mf{g}}_{-})\otimes A_{\lambda} 
 \cong V_{\lambda, k}\otimes V_{\lambda^{*}} $ and 
 $ U(\bar{\hat{\mf{g}}}, \bar{k}) \otimes_{U(\bar{\hat{\mf{g}}}_{\geq 0}, \bar k)} A_{\lambda} 
 = U(\bar{\hat{\mf{g}}}_{-})\otimes A_{\lambda} 
 \cong V_{\lambda}\otimes V_{\lambda^{*}, \bar{k}}$. 

\begin{lemma}\label{B1}
$B$ can be equipped with a $\hat{\mf{g}}$- (resp. $\bar{\hat{\mf{g}}}$-)module structure 
denoted by $\rho$ (resp. $\tilde{\rho}$), 
such that 
$\rho|_{B_{\lambda}} \cong (V_{\lambda^{*}, k}\otimes V_{\lambda})^{c}$ 
(resp. $\tilde{\rho}|_{\widetilde{B_{\lambda}}} \cong (V_{\lambda^{*}}\otimes V_{\lambda, \bar{k}})^{c}$). 
\end{lemma}

\begin{proof}
$(,)$ induces a $\hat{\mf{g}}$- (resp. $\bar{\hat{\mf{g}}}$-)homomorphism 
$\phi$ (resp. $\psi$) from $U\mathcal{A}$ to 
$(\oplus_{\lambda} V_{\lambda, k}\otimes V_{\lambda^{*}})^{c}$ 
(resp. $(\oplus_{\lambda} V_{\lambda}\otimes V_{\lambda^{*}, \bar{k}})^{c}$). 
Let $J$ (resp. $\bar{J}$) be the augmentation ideal in $U(\hat{\mf{g}}_{-})$ 
(resp. $U(\bar{\hat{\mf{g}}}_{-})$), 
then $U\mathcal{A} = (\bar{J} \otimes B)\oplus B$ 
(resp. $U\mathcal{A} =(J \otimes B)\oplus B$) 
and it is obvious that 
$\phi |_{\bar{J} \otimes B} =0$ (resp. $\psi |_{J \otimes B} =0$). 
In addition by the proof of Proposition \ref{contragrediant}, 
the pairing between $U(\hat{\mf{g}}_{-}) \otimes A_{\lambda}$ and $B_{\lambda^{*}}$ 
(resp. between $U(\bar{\hat{\mf{g}}}_{-})\otimes A_{\mu}$ and $\widetilde{B_{\mu^{*}}}$) 
is non-degenerate, 
hence the induced maps 
$\phi|_{B_{\lambda^*}}: B_{\lambda^{*}} \to (V_{\lambda, k}\otimes V_{\lambda^{*}})^{c}$ 
(resp. $\psi|_{\widetilde{B_{\mu^*}}}: \widetilde{B_{\mu^*}} \to  (V_{\mu}\otimes V_{\mu^{*}, \bar{k}})^{c}$) are isomorphisms of vector spaces 
(since all $B_{\lambda^*}$'s and  $\widetilde{B_{\mu^*}}$'s 
are finite dimensional on each level). 
Hence  
$\phi |_{B}: B \to \oplus_{\lambda} (V_{\lambda, k}\otimes V_{\lambda^{*}})^{c}$ 
(resp. $\psi |_{B}: B \to \oplus_{\mu} (V_{\mu}\otimes V_{\mu^{*}, \bar{k}})^{c}$) 
is an isomorphism of vector spaces. 
Therefore the $\hat{\mf{g}}$- (resp. $\bar{\hat{\mf{g}}}$-)module structure of the right hand side can be carried over to the left.  
\end{proof}

\begin{remark}
$B$ is not a $\hat{\mf g}$- (or $\bar { \hat {\mf g}}$-)submodule of $U\mathcal A$, 
but only closed with repsect to the action of 
$\hat {\mf g}_{\geq 0}$ (or $\bar {\hat {\mf g}}_{\geq 0}$). 
The previous lemma means that 
$x(m) b \equiv \rho(x(m)) b \, (\text{mod  } \bar{J} \otimes B)$ and 
$\bar{y}(n) b \equiv \tilde{\rho}(\bar{y}(n)) b \, (\text{mod  } J \otimes B)$
for any $x(m) \in \hat {\mf g}$, $\bar y(n) \in \bar {\hat {\mf g}}$ and $b\in B$. 
\end{remark}

\begin{lemma}\label{g}
Suppose $k, \bar{k}\notin \mathbb{Q}$, then
\begin{enumerate}
\item
for any $\lambda \in P^{+}$ the pairing 
$(, ): (U(\hat{\mf{g}}_{-})\otimes A_{\lambda}) \times (U(\hat{\mf{g}}_{-})\otimes A_{\lambda^{*}}) \to \mathbb{C}$ (resp. $(, ): (U(\bar{\hat{\mf{g}}}_{-})\otimes A_{\lambda}) \times (U(\bar{\hat{\mf{g}}}_{-})\otimes A_{\lambda^{*}})\to \mathbb{C}$) is nondegenerate.
\item 
$V_{\lambda, k} \otimes V_{\lambda^{*}, \bar{k}}$ is irreducible as a $\hat{\mf{g}} \oplus \bar{\hat{\mf{g}}}$-module. 
\end{enumerate}
\end{lemma}

\begin{proof}
Obvious since $V_{\lambda, k}$ and $V_{\lambda^{*}, \bar{k}}$ 
are irreducible when $k, \bar{k}\notin \mathbb{Q}$. 
\end{proof}

\begin{remark}
For generic values of $k$ and $\bar{k}$, 
by (2) of the previous lemma, 
the $\hat{\mf{g}} \oplus \bar{\hat{\mf{g}}}$-submodule 
$U\mathcal{A}^{s}$ generated by $A$ 
is as big as 
$U(\hat{\mf{g}}_{-}) \otimes U(\bar{\hat{\mf{g}}}_{-}) \otimes A$. 
It is not difficult to show that it is indeed the whole $U\mathcal{A}$ 
by comparison of the characters. 
But here we prove it more directly for the reason that the computations we use will be helpful later for the non-generic situation. 
\end{remark}

\begin{lemma}\label{lem2}
Let $b\in B$, then 
\begin{enumerate}
\item 
if $m \geq 0$, then $\tau_i(m) b\in B$ (resp. $\bar{\tau_i}(m) b\in B$); 
\item 
if $m >0$, then 
$\tau_{i}(-m) b = \sum_{j, 1\leq n\leq m} \bar{\tau_j}(-n) b^{j, n} + \rho( \tau_{i}(-m)) b$ 
for some $b^{j, n}\in B$ 
(resp. 
$\bar{\tau_{i}}(-m) b = \sum_{j, 1\leq n\leq m} \tau_{j}(-n) \widetilde{b^{j, n}} + \tilde{\rho} (\bar{\tau_{i}}(-m)) b$ 
for some $\widetilde{b^{j,n}} \in B$);  
\item 
generally for any $m >0$ and $n \geq 0$, one has 
$U(\hat{\mf{g}}_{-}) [m] \otimes B_n \subset 
\sum_{i=1}^{m} U(\bar{\hat{\mf{g}}}_{-})[i] \otimes B_{n+m -i} + B_{m+n}$ 
(resp. $U(\bar{\hat{\mf g}}_-) [m] \otimes B_n \subset 
\sum_{i=1}^{m} U(\hat{\mf g}_-)[i] \otimes B_{n+m -i} +  B_{m+n}$); 
\item 
define a map $t: U(\hat{\mf{g}}_{-})\otimes A \to U(\bar{\hat{\mf{g}}}_{-})\otimes A$ 
as follows: 
$a\mapsto a$; 
$\tau_{i_{1}} (-m_{1})\cdots \tau_{i_{l}}(-m_{l}) a \mapsto \sum_{j_{1}, \cdots, j_{l}} 
\bar{\tau_{j_{l}}} (-m_{l})\cdots \bar{\tau_{j_{1}}} (-m_{1}) ( a^{j_{l}i_{l}} \cdots a^{j_{1}i_{1}} a)$ 
for any $a\in A$, $m_{1},  \ldots, m_{l} >0$, 
then $t$ is bijective and $( \cdot, b) = (t (\cdot) , b)$ for any $b\in B$. 
\end{enumerate}
\end{lemma}

\begin{proof}
(1) is already proved in Proposition \ref{PBW}.
(2) follows from (8) and (9) of Lemma \ref{ca}. 
In fact we have 
$\tau_{i}(-m) b 
= \sum_{j} \{\bar{\tau_{j}} (-1) (b \partial^{(m-1)} a^{ji}) 
+ \bar{\tau_{j}}(-2)( b\partial^{(m-2)} a^{ji} )
+ \cdots + \bar{\tau_{j}}(-m) (b a^{ji})\} 
+ \rho( \tau_{i}(-m)) b$ 
(resp. 
$\bar{\tau_{i}}(-m) b 
= \sum_{j} \{\tau_{j}(-1) (b\partial^{(m-1)} a^{ij}) 
+\tau_{j}(-2) (b\partial^{(m-2)} a^{ij}) 
+\cdots + \tau_{j}(-m) (b a^{ij}) \} 
+ \tilde{\rho} (\bar{\tau_{i}}(-m)) b$). 
Set $b^{j, n} = b \partial^{(m-n)} a^{ji} $ 
(resp. $\widetilde{b^{j,n}} = b\partial ^{(m-n)} a^{ij}$), 
then for homogeneous $b$, 
$b^{j, n}\in B_{m- n + \text{deg} b}$ 
(resp. $\widetilde{b^{j, n}} \in B_{m -n +\text{deg} b}$).
(3) can be proved by induction. 
(4) follows from (3) and the fact that 
$(U(\bar{\hat{\mf{g}}}_{-}) \otimes (\oplus_{n\geq 1} B_{n}), B) = 0$.
\end{proof}

\begin{flushleft}
\bf{Proof of Theorem \ref{main2}}:
\end{flushleft}

To prove that $U\mathcal{A} = U\mathcal{A}^{s} $ when $k, \bar k \notin \mathbb Q$, 
it suffices to show that $B \subseteq U\mathcal{A}^{s}$. 
Clearly $B_{0} = A\subseteq U\mathcal{A}^{s}$. 
Suppose $B_{k} \subseteq U\mathcal{A}^{s}$ for all $k\leq n$, 
we want to show that $B_{n+1}\subseteq U\mathcal{A}^{s}$. 
For any  $b_{n+1}\in B_{n+1}$, by Lemma \ref{B1} and \ref{g}, 
there exists a $v\in (U(\hat{\mf{g}}_{-})\otimes A)_{n+1}$ 
such that $(b_{n+1} - v, U(\hat{\mf{g}}_{-})\otimes A) =0 $, 
i.e. $b_{n+1}-v \in \bar{J}\otimes B$.  
By Lemma \ref{lem2}, only elements from $B_{\leq n}$ appear in $b_{n+1} -v\in \bar{J}\otimes B$, 
hence by assumption $b_{n+1} -v \in U\mathcal{A}^{s}$, 
therefore $b_{n+1}\in U\mathcal{A}^{s}$. 

To prove that $U\mathcal{A}\cong V$ as vertex operator algebras, 
it suffices to show that there is only one possible VOA structure on the space 
$\oplus_{\lambda} V_{\lambda, k}\otimes V_{\lambda^{*}, \bar{k}}$ ($k \notin \mathbb{Q}$)
if we require that $V_{0, k}\otimes V_{0, \bar k}$
yields the $\hat{\mf{g}}\oplus \bar{\hat{\mf{g}}}$-module structure of it 
as a vertex subalgebra, 
and the restriction of the  $_{(-1)}$ operation to 
$\oplus_{\lambda} V_{\lambda}\otimes V_{\lambda}^{*} \cong \mathcal{R}(G)$
is multiplication of regular functions. 
Using the nondegenerate symmetric bilinear form on  
$\oplus_{\lambda} V_{\lambda, k}\otimes V_{\lambda^{*}, \bar{k}}$ ($k \notin \mathbb{Q}$)
determined by $x(n)^{*}= - x(-n)$, $\bar{y}(m)^{*}= - \bar{y}(-m)$ and Definition \ref{bfonA}, 
it is not difficult to see that $u_{(n)}v$ is uniquely determined by the requirements for any 
$u, v\in \oplus_{\lambda} V_{\lambda, k}\otimes V_{\lambda^{*}, \bar{k}}$ and $n\in \mathbb{Z}$. 
Thus Theorem \ref{main2} is proved.

\begin{flushright}
$\square$
\end{flushright}

For the rest of Section 3 we will show that the Zhu's algebra for $U\A_{\mf{g}, k}$ is isomorphic to the algebra of differential operators on the Lie group $G$.  

In [Zh], Y. Zhu defined an associative algebra $A(V)$ for any vertex operator algebra $V$ and proved the one-to-one correspondence between irreducible $A(V)$-modules and irreducible $V$-modules. 
As a vector space $A(V)$ is a quotient space of $V$. 
More specifically let $O(V)$ be the linear span of elements of the form 
$\text{Res}_{z}(Y(u, z) \frac{(z+1)^{\text{wt} u}}{z^{2}} v)$
where $u, v\in V$ with $u$ homogeneous, and set $A(V)=V/O(V)$. 
Define a multiplication in $V$ as follows: 
$$ u  \ast v = \text{Res}_{z} (Y(u, z) \frac{(z+1)^{\text{wt} u}}{z} v), $$
then it induces the (associative) multiplication on the quotient $A(V)=  V/O(V)$ 
and the image of the vacuum $\mathbf{1}$ in $A(V)$ becomes the identity element. 
For example, the Zhu's algebra for the vertex operator algebra $V_{0, k}$ is isomorphic to the universal enveloping algebra of $\mf g$ ([FZ]). 

Denote by $[v]\in A(V)$ the image of $v\in V$. 

\begin{lemma}\label{tiny}
For any $a\in V_{0}$, $x\in V_{1}$, $v\in V$ and $n\geq 0$, one has
\begin{enumerate}
\item 
$[a] \ast [v] =  [a_{(-1)} v]$, $[x] \ast [v] = [ (x_{(-1)} + x_{(0)}) v]$ and $[v] \ast [x] = [ x_{(-1)} v]$; 
\item 
$a_{(-2-n)} v, (x_{(-1-n)} + x_{(-2-n)}) v \in O(V)$.  
\end{enumerate}
\end{lemma}

\begin{proof}
(1) follows from the definition and Lemma 2.1.7 in [Zh]. 
(2) follows from Lemma 2.1.5 in [Zh].
\end{proof}

\begin{lemma}\label{small}
For the vertex operator algebra $U\A = U\A_{\mf{g}, k}$ one has 
\begin{enumerate}
\item 
$U(\hat{\mf{g}}_{-}) \otimes B_{\geq 1} \subset O(U\A)$; 
\item
$[x_{1}(-n_1-1) \cdots x_{l}(-n_l-1) a] = (-1)^{n_{1}+\cdots n_{l}} [a] \ast [x_{l}(-1) \mathbf{1}] \ast \cdots \ast [x_{1}(-1) \mathbf{1}]$ 
for any $x_{i} \in \mf{g}$, $n_{i}\geq 0$ and $a\in A$. 
\end{enumerate}
\end{lemma}

\begin{proof}
It follows from Lemma \ref{AB} and (2) of Lemma \ref{tiny} that $[\omega_{r}] =0$ for any $r$. 
Moreover for any $n \geq 0$, $v\in U\A$,  by Lemma \ref{tiny} we have 
$[{\omega_{r}}_{(-n-1)} v] = (-1)^{n} [{\omega_{r}}_{(-1)} v] = (-1)^{n} [v] \ast [\omega_{r}] = 0$, 
hence $[b] = 0$ for any $b\in B_{\geq 1}$. 
Again by (2) of Lemma \ref{tiny}, we have 
$[x_{1}(-n_{1}-1) \cdots x_{l}(-n_{l}-1) b] 
= (-1)^{n_{1}} [ x_{1}(-1)\cdots x_{l} (-n_{l}-1) b]
=(-1)^{n_{1}} [x_{2}(-n_{2}-1)\cdots x_{l}(-n_{l}-1) b] \ast [x_1(-1) \mathbf{1}] 
= \cdots
= (-1)^{n_{1}+\cdots n_{l}} [b] \ast [x_{l}(-1)\mathbf{1}]\ast \cdots \ast [x_{1}(-1)\mathbf{1}] $
for any $x_i\in \mf g$, $n_i\geq 0$ and $b\in B$. 
The lemma now follows. 
\end{proof}

Let $\mf{v} = A \oplus \mf{g} $ be the Lie algebra where the Lie bracket is given  by 
$$[ a + x, a' + y ] = (x \cdot a' - y \cdot a) + [x, y]$$ 
for any $a, a'\in A$, $x, y\in \mf{g}$.  
Since $A$ is an abelian Lie subalgebra of $\mf{v}$,  
we have the decomposition 
$U(\mf{v}) = \textrm{Sym}_{\C} A \otimes U(\mf{g})$. 
Let $I$ be the ideal of $U(\mf{v})$ generated by elements of the form 
$ a \cdot a' - aa'$ where $a \cdot a'$ denotes the formal multiplication in $\textrm{Sym}_{\C} A$ while $aa'$ denotes the actual multiplication of functions. 
Let $\mathcal{D} = U(\mf{v})/ I$, 
then $\mathcal{D} \cong A \otimes U(\mf{g})$ as a vector space where the isomorphism is induced by multiplication. In fact $\mathcal D$ is the algebra of differential operators on the Lie group $G$. 

\begin{prop}
$A(U\A) \cong  \mathcal{D}$. 
\end{prop}

\begin{proof}
By Lemma \ref{small}, $A(U\A)$ is generated as an associative algebra  by elements of the type 
$[a], [x(-1)\mathbf{1}]$ where $a\in A, x\in \mf{g}$. 
By definition $[a] \ast [a'] = [a_{(-1)} a'] = [aa'] = [a'] \ast [a]$. 
Moreover by Lemma 2.1.7 in [Zh], we have 
$x(-1)\mathbf{1} \ast a - a \ast x(-1)\mathbf{1}\equiv \textrm{Res}_{z} (Y(x(-1)\mathbf{1}, z) a ) \equiv x(0)a \equiv x\cdot a \,\, \textrm{mod  } O(U\A)$, i.e. $[x(-1)\mathbf{1}] \ast [a] - [a] \ast [x(-1)\mathbf{1}] = [x \cdot a]$. 
Hence there exists a surjective homomorphism $\mathcal{D} \twoheadrightarrow A(U\A)$.
To see that it is also injective it suffices to show that $O(U\A)$ is the linear span of $U(\hat{\mf{g}}_{-}) \otimes B_{\geq 1}$ and elements of the type $(x(-n-2) + x(-n-1)) v$ 
where $x\in \mf{g}$, $n\geq 0$ and  $v\in U\A$, 
which can be proved by induction. 
\end{proof}

\section{integral central charges}
Finally we would like to discuss the case when $k,\bar{k}\in \mathbb{Z}$. 
This case is more subtle and undoubtedly more interesting (cf. [FS] 4.3).  
The Fock space realization in [FS] exists for any $\varkappa\neq 0$, i.e. excluding the critical levels $k=\bar{k}=-h^{\vee}$, and 
the construction of the vertex envelop of the vertex algebroid $\mathcal A_{\mf g, k}$ studied in Section 3 works for any $k\in \mathbb{C}$.
When $k\in \mathbb Z$, the vertex operator algebra $U\mathcal A$ is the same size as in the generic case, but the $\hat{\mf g} \oplus \bar{\hat {\mf g}}$-module structure is much more complicated. 
Instead of having a nice decomposition, it has linkings determined by the shifted action of the affine Weyl group. We have the following results based on the work done in Section 3.

The Weyl modules with rational central charges are usually reducible. We say that a vector $v\in V_{\lambda, k}$ is singular if $\hat{\mf{g}}_{>0} \cdot v=0$. The singular vectors are closed under the $\mf g$-action. 

Let $S=(U(\hat{\mf{g}}_{-})\otimes A)\cap (U(\bar{\hat{\mf{g}}}_{-})\otimes A)$. 

\begin{lemma}
$S = U\mathcal{A}^{\hat{\mf{g}}_{>0}\oplus \bar{\hat{\mf{g}}}_{>0}} = \{ v\in U(\hat{\mf{g}}_{-})\otimes A: t (v) =v \}$. 
\end{lemma}

\begin{proof}
The first equality follows from Proposition \ref{PBW}.
The second one follows from (4) of  Lemma \ref{lem2} because the map $t$ is uniquely determined by the pairing property. 
\end{proof}

$S$ is $\mathbb{N}$-graded and $S_{0} = A$. 
If $k, \bar{k}\notin \mathbb{Q}$, then $S=A$, otherwise $S$ is much bigger and contains all the singular vectors. 

\begin{prop} \label{4.2}
If $X \subset V_{\lambda, k}$ is a subspace of singular vectors and $X\cong V_{\mu}$ as a $\mf{g}$-module, then there exists a singular subspace $Y\subset V_{\mu^{*}, \bar{k}}$ such that
$Y\cong V_{\lambda}^{*}$ as a $\mf{g}$-module 
and $X\otimes V_{\lambda}^{*} = V_{\mu} \otimes Y$ in $U\mathcal{A}$. 
\end{prop}

\begin{proof}
Suppose that $X \subset V_{\lambda, k}$ is a singular subspace and $X\cong V_{\mu}$ as a $\mf{g}$-module, then $X \otimes V_{\lambda}^{*} \subset (U\mathcal{A})^{\hat{\mf{g}}_{>0}} = U(\bar{\hat{\mf{g}}}_{-})\otimes A$. It follows that $X \otimes V_{\lambda}^{*} \subset U(\bar{\hat{\mf{g}}}_{-})\otimes A_{\mu}$, hence $X\otimes V_{\lambda}^{*} = V_{\mu}\otimes Y$ for some singular subspace $Y\subset V_{\mu^{*}, \bar{k}}$ and $Y\cong V_{\lambda^{*}}$ as a $\mf{g}$-module. 
\end{proof}

Let $\lambda, \mu\in P^{+}$ be as in Proposition \ref{4.2}, then the $\hat{\mf{g}}\oplus \bar{\hat{\mf{g}}}$-submodules generated by $V_{\lambda}\otimes V_{\lambda}^{*}$ and $V_{\mu}\otimes V_{\mu}^{*}$ have a nontrivial intersection, 
therefore $U\mathcal{A}^{s}$ (the $\hat{\mf{g}}\oplus \bar{\hat{\mf{g}}}$-submodule generated by $A$) 
is strictly smaller than $U\mathcal{A}$.

In view of the equivalence of categories of representations of the affine Lie algebra and representations of the quantum group (cf. [KL1-4]), 
and the Arkhipov functor which transforms modules with a Weyl filtration of level $k$ into modules with a Weyl filtration of dual level $\bar k$ (cf. [A], [So]), 
we conjecture that under equivalence of categories, 
the $\hat{\mf{g}}\oplus \bar{\hat{\mf{g}}}$-module structure of $U\mathcal{A}$ with noncritical integral central charges is the same as the bimodule structure of the regular representation of the corresponding quantum group at roots of unity (see [Z]). 
For $k\geq 0$, $\bar{k}\leq -2 h^{\vee}$, we conjecture that 
$U\mathcal{A}$ admits an increasing filtration of $\hat{\mf{g}}\oplus \bar{\hat{\mf{g}}}$-submodules having $V_{\lambda, k}\otimes V_{\lambda, \bar{k}}^{c}$ as successive quotients and a decreasing filtration with $V_{\lambda, k} ^{c}\otimes V_{\lambda, \bar{k}}$ as successive quotients, where $V_{\lambda, k}^{c}$ is the contragredient module of $V_{\lambda, k}$ defined by the antiautomorphism $x(n)\mapsto -x(-n)$, $\underline c\mapsto \underline c$ of $\hat{\mf{g}}$.

\end{document}